%% file: main.tex
\documentclass[article]{siamart}
\usepackage{amsmath,amssymb}
\usepackage{caption}
\usepackage{hyperref}
\usepackage[utf8]{inputenc}
\usepackage{mathtools}
\usepackage{color,graphicx}
\usepackage{thmtools, thm-restate}
\usepackage{enumitem}
\usepackage{algorithmicx}
\usepackage{algpseudocode}
\usepackage{tikz}
\usetikzlibrary{patterns}
\usetikzlibrary{arrows}
\usetikzlibrary{3d,calc}

\input{macros.tex}

\newcommand{\TheTitle}{Model reduction of a parametrized scalar hyperbolic 
                    conservation law using displacement interpolation}

\newcommand{\remove}[1]{\textcolor{red}{[Text removed.]}}%
\graphicspath{{./figures/}{./slides/figures/}{./plots/savetikz/}}
\begin{document}

\ifpdf
\DeclareGraphicsExtensions{.pdf, .jpg, .tif}
\else
\DeclareGraphicsExtensions{.eps, .jpg}
\fi

\title{\TheTitle}

\author{Donsub Rim%
  \thanks{Department of Applied Physics and Applied Mathematics,  %
  Columbia University, New York, %
  NY 10027 (\email{dr2965@columbia.edu}, \email{kyle.mandli@columbia.edu}).}%
  \and %
  Kyle T. Mandli\footnotemark[1]  %
}
\maketitle

\begin{abstract} 
We propose a model reduction technique for parametrized partial differential 
equations arising from scalar hyperbolic conservation laws.
The key idea of the technique is to construct basis functions that are local in 
parameter and time space via displacement interpolation. The construction
is motivated by the observation that the derivative of solutions to 
hyperbolic conservation laws satisfy a contractive property with respect to the 
Wasserstein metric 
[Bolley et al. \emph{J. Hyperbolic Differ. Equ.} 02 (2005), pp. 91-107]. 
We will discuss the approximation properties of the displacement
interpolation, and show that it can naturally complement linear interpolation. 
Numerical experiments illustrate that we can successfully achieve the model 
reduction of a parametrized Burgers' equation, and that the reduced order model 
is suitable for performing typical tasks in uncertainty quantification.
\end{abstract}

\section{Introduction}

Although reduced order models (ROMs) have been successfully constructed for 
various partial differential equations (PDEs), existing projection-based methods 
\cite{berkooz93,holmes96,benner15} are often ineffective when applied 
directly to problems governed by hyperbolic PDEs, and the development of ROMs 
for hyperbolic conservation laws still is in its early stages.
The main obstacle stems from the fact that the energy of the solution is
typically concentrated at different spatial locations for different time or 
parameter values \cite{rowley00,amsallem,carlberg15}. Put in other terms, 
the Kolmogorov $N$-width decays slowly with respect to the dimension of the 
reduced basis for hyperbolic problems \cite{crb-book,gpc-book10}. 
Due to this slow decay, the construction of a global basis appears difficult. 

In this paper, we build on previous works \cite{rim17reversal, rim17mr}
to propose a procedure which constructs a basis that is \emph{local} 
in the parameter-time space, that neither requires additional queries to the 
high-fidelity model (HFM), nor utilizes adaptive procedures during the online 
stage of the ROM. In our approach, 
rather than finding a low-rank approximation for the snapshots directly,
we compute the \emph{transport maps} between the snapshots in order to 
find a low-dimensional structure in these transport maps.
Our technique shares common features with recently proposed methods for 
extracting low-dimensional transport structures to be used for model reduction 
\cite{schulze15,schulze18,rim17reversal,welper17p,welper17,rim17mr}.

This key idea is closely related to problems arising 
in optimal transport. In particular, we will make use of an interpolation 
procedure called \emph{displacement interpolation} in the 
optimal transport literature \cite{villani2008,villani03}.
In a single spatial dimension (1D), this transport map can be computed
explicitly by employing the simplest solution to the Monge-Kantorovich 
problem \cite{Kantorovich} called \emph{monotone rearrangement} \cite{Benamou}.
We will derive a variant of monotone rearrangement, 
guided by a crucial relation between scalar conservation laws and optimal 
transport observed in \cite{bolley05}.

In justifying our construction, we find that the displacement interpolant
itself has very general approximation properties. In some sense, this
approximation is \emph{dual} to linear approximation, when one views it as a 
linear approximation in the dependent variable in a natural coordinate transform.
We make this explicit in \cref{sec:approx} and propose a more general form
for approximation, one that expands the transport map in some suitable
basis \cref{eq:x_sum}, as opposed to the usual linear approach \cref{eq:u_sum}. 
The latter underlies many prevailing numerical methods in dimensionality
reduction, e.g., singular value decomposition (SVD) \cite{horn12,golub96},
basis pursuit \cite{chen01} or  generalized polynomial chaos (gPC) \cite{gpc-book10}.

This paper is organized to guide the reader through the step-by-step 
construction of the reduced order model. The key steps are as follows.
\begin{center}
\begin{tabular}{rl}
Step 1.& Discretization of the parameter-time space 
(\cref{sec:def})  \\
Step 2.& Construction of local basis using displacement interpolation 
(\cref{sec:basis}) \\
Step 3.& Construction of the ROM via projection 
 (\cref{sec:rom})  \\
Step 4.& Dimensionality reduction using standard techniques 
        (\cref{sec:rom})  \\
\end{tabular}
\end{center}
In summary, we first obtain a discretization of the parameter-time space that 
satisfies the so-called \emph{signature condition} (\cref{cond:signature}),
then we construct a basis that is local with respect to this discretization.
At this stage, we may construct a ROM by projecting the PDE to this basis.
Finally, we further reduce this basis using standard projection techniques
such as proper orthogonal decomposition (POD).

In \cref{sec:approx} we will discuss the approximation properties of the 
displacement interpolation that will yield more insight into the local basis
construction procedure. Then in \cref{sec:uq} we will illustrate that the ROM 
can then be used reliably for common tasks
in uncertainty quantification (UQ) by applying Monte Carlo sampling
on the constructed ROM.

\section{Parametrized Burgers' equation} \label{sec:par_burg}

Throughout this paper, we will consider the example of a parametrized Burgers'
equation, which was used as a benchmark problem in 
\cite{carlberg13,carlberg15,choi17,reweinski98,rewienski06}: 
\begin{equation}
\begin{aligned}
\frac{\partial u}{\partial t} + \frac{\partial}{\partial x}\left( \frac{1}{2}
u^2 \right) &= 0.02 e^{ -\mu_2 x }, 
\quad \text{ for } (x,t) \in (0,100) \times (0,50)\\
u(x,0) &= 0,\\
u(0,t) &= \mu_1,\\
\end{aligned}
\label{eq:par_burg}
\end{equation}
where $\mu_1$ is the incoming boundary condition.
The values $\mu_1$ and $\mu_2$ are drawn from the parameter
space $\cM$,
\beq
\bmu = (\mu_1,\mu_2) \in \cM = [3,9] \times [0.02,0.075].
\label{eq:cM}
\eeq
We will denote the solution to \cref{eq:par_burg} with given parameter values
$\bmu$ by
\beq
u(x,t; \bmu) = u(x,t; \mu_1,\mu_2).
\eeq
The HFM will use the finite volue method (FVM) \cite{fvmbook} to solve
the equation, with the uniform finite volume cell width $\Delta x = 0.4$ 
amounting to $N=250$ total number of cells, and time-step $\Delta t = 0.0125$.
The finite volume cell at time $t_n$ will be denoted by,
\beq
U_{i,n} := \frac{1}{\Delta x} \int_{x_{i-1/2}}^{x_{i+1/2}} u(x,t_n) \dx.
\label{eq:fvcell}
\eeq
The time-step update for the HFM is given by
\beq
    U_{i,n+1} = U_{i,n} - \frac{\Delta t}{\Delta x}
        (\cF[U_{i+1,n}] - \cF[U_{i,n}])
                + \Delta t \, \left( 0.02 e^{\mu_2 x_i} \right),
\label{eq:hfm_update}
\eeq
where Godunov's method will be used to compute the numerical flux $\cF$. 
Since our solution will not have positive jump discontinuities, we will not 
apply the entropy fix in the time-step updates, even though it would be
easy to include the fix if needed.

\begin{figure}
    \centering
        \begin{tabular}{cc}
        {\small $(\mu_1,\mu_2) = (9,0.02)$}& 
        {\small $(\mu_1,\mu_2) = (3,0.075)$}\\
        \includegraphics[width=0.45\textwidth]{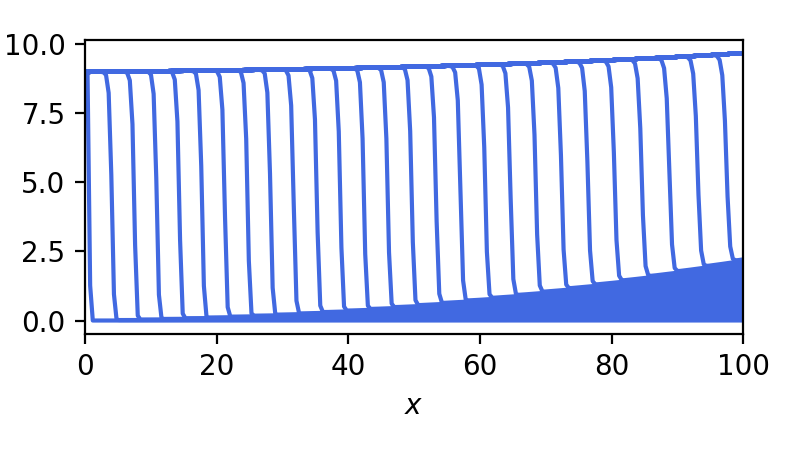} &
        \includegraphics[width=0.45\textwidth]{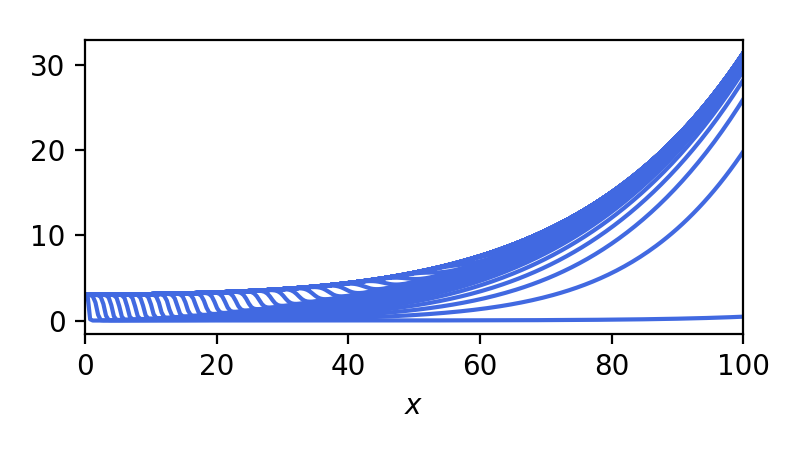} \\
        {\small $(\mu_1,\mu_2) = (3,0.02)$}& 
        {\small $(\mu_1,\mu_2) = (6,0.05)$}\\
        \includegraphics[width=0.45\textwidth]{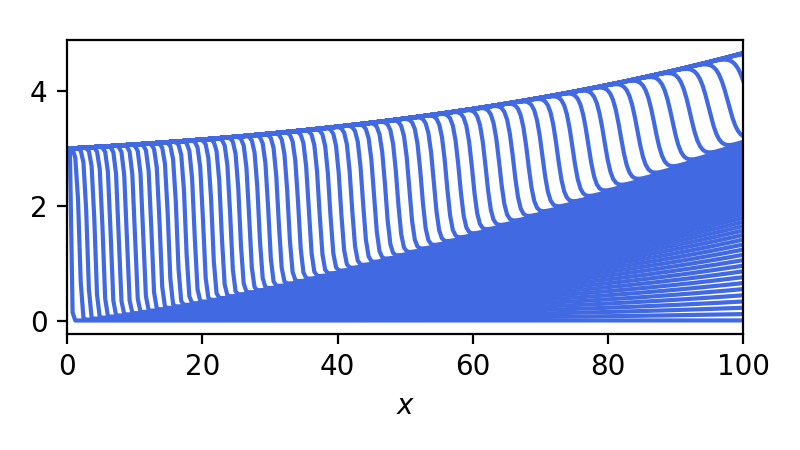} &
        \includegraphics[width=0.45\textwidth]{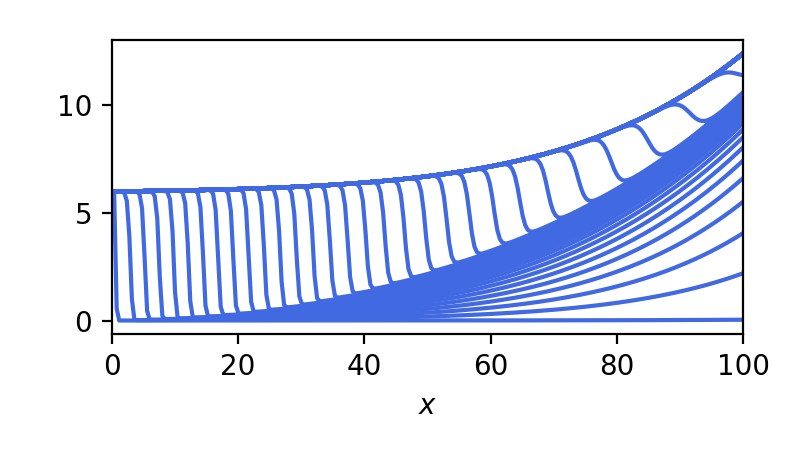} \\
        \end{tabular}
    \caption{HFM solutions using the finite volume method \cref{eq:hfm_update}
             for various parameter values, plotted every 60 time-steps.}
    \label{fig:hfm_sols}
\end{figure}

The HFM will be run for 9 parameter values in the set $\cM_0$,
\beq
(\mu_1,\mu_2) \in \cM_0 \equiv \{3,6,9\} \times \{0.02, 0.05, 0.075\}.
\label{eq:M0}
\eeq
A few of the computed solutions are plotted in Figure \ref{fig:hfm_sols}.

\section{Definitions and notations}\label{sec:def}

In this secton, we will set up definitions and notations related to the
discretization of the parameter-time space, as well as those related to 
displacement interpolation. They will be used throughout this paper.

\subsection{Parameter-time elements}

Since our ROM will depend on basis functions that are local in
parametric variables $\{\mu_1,\mu_1\}$ and time variable $t$, we will first discretize the parameter-time space. We will 
make use of a Delaunay triangulation for its simplicity, 
although this is not the only choice. 

Let us first denote the entire parameter-time space by $\cM_T$,
\beq
\begin{aligned}
\cM_T &= \left\{(\mu_1,\mu_2,t): \mu_1 \in [3,9], 
                                 \mu_2 \in [0.02,0.075], 
                                     t \in [0, 50] \right\}. 
\end{aligned}
\eeq
We will proceed by partitioning $\cM_T$ as follows.
\begin{itemize}
    \item We compute the Delaunay triangulation over the parameter space 
    $\cM$ \cref{eq:cM}.
    The $\ell$-th triangle in this triangulation will be denoted by $\cT^\ell$ 
    and is designated by three points 
    $\{\bmu_{j_{\ell,1}},\bmu_{j_{\ell,2}},\bmu_{j_{\ell,3}}\} \subset \cM_0$
    in the parameter space. That is,
    \beq
        \cT^\ell :=
            (\text{convex hull of }
                    \{\bmu_{j_{\ell,1}},\bmu_{j_{\ell,2}},\bmu_{j_{\ell,3}} \}
            )
    \label{eq:tri}
    \eeq
    The triangulation we will use is plotted in \cref{fig:tri}.
    \ezfigure{width=0.6\textwidth}{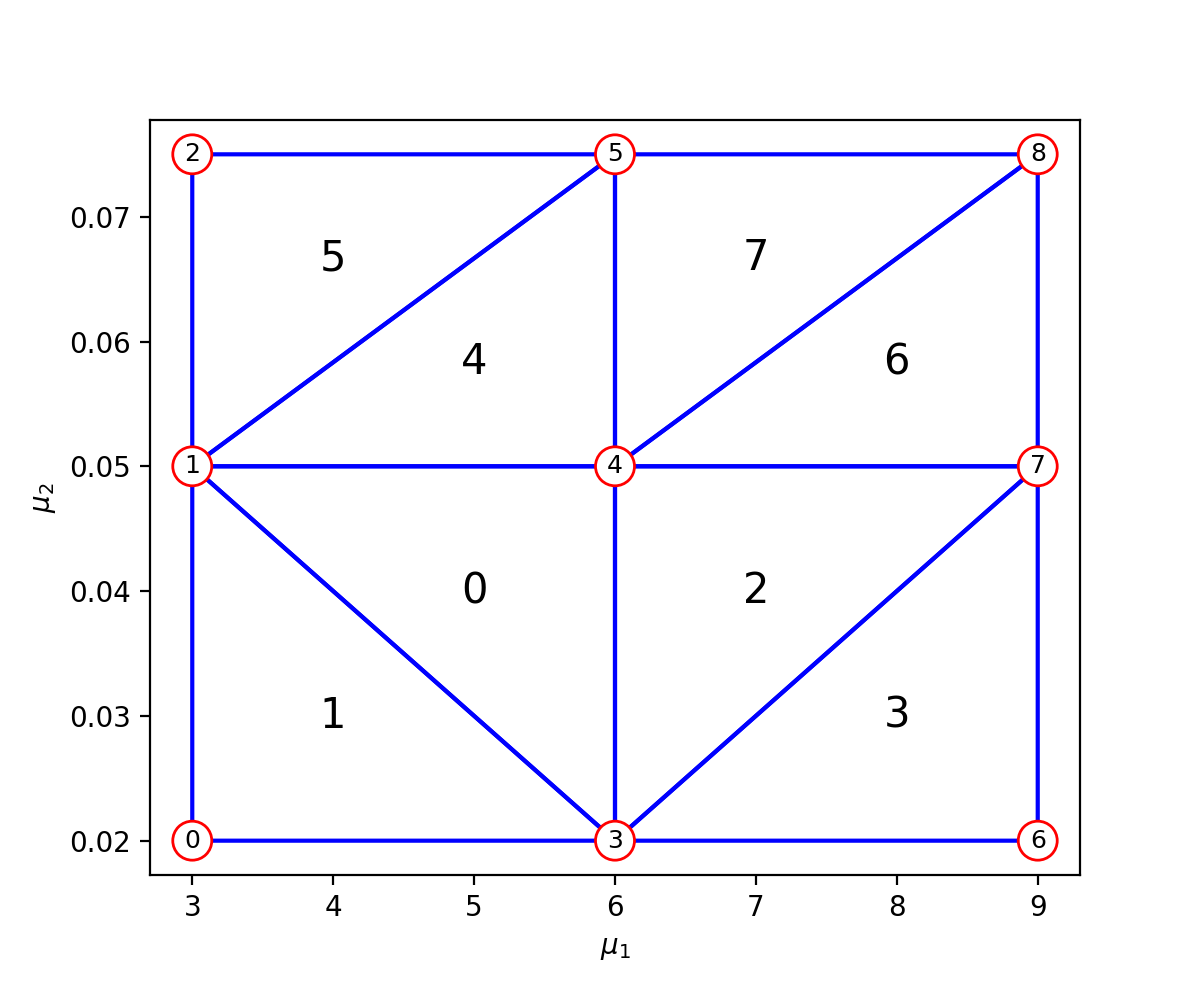}
             {Plot of the triangulation defined in \cref{eq:tri}
              used to discretize the parameter space $\cM$ \cref{eq:cM}.
              The index $\ell$ is displayed in the interior of each $\cT^\ell$,
              and the index of parameter values in $\cM_0$ \cref{eq:M0}
              is displayed at the nodes of these triangles.
             }%
             {fig:tri}
    \item We will extend the triangle $\cT^\ell$ as an element over the 
          parameter-time space by selecting a subset of the time-steps 
          $\{t_{n_m}\} \subset \{t_n\}$
    \beq
        \cE^\ell_m
        := \cT^\ell \times [t_{n_m},t_{n_{m+1}}).
    \label{eq:cE}
    \eeq
    A diagram of such an element is plotted in \cref{fig:elem}.
    \item Nodes (or vertices) of $\cE^\ell_m$ are the points in parameter-time
          space,
    \beq
        \cN(\cE^\ell_m) := \{
                \bmu_{j_{\ell,1}},\bmu_{j_{\ell,2}},\bmu_{j_{\ell,3}} 
                           \}
                \times \{t_{n_m}, t_{n_{m+1}}\}.
        \label{eq:nodes}
    \eeq
    Note that each of these nodes can be related to the solution to
    \cref{eq:par_burg}. 
    To each member of $\cN(\cE^\ell_m)$ there corresponds a function of the 
    spatial variable $x$,
    \beq
        u(x,t_*; \bmu_*) \where (\bmu_*, t_*) \in \cN(\cE^\ell_m).
        \label{eq:uxtmuN}
    \eeq
    Since the HFM solution will be computed for each 
    $\bmu \in \cM_0$ \cref{eq:M0} during the offline-phase, we assume that
    we have an accurate approximation to the $x$-dependent function 
    \cref{eq:uxtmuN} via the HFM \cref{eq:hfm_update}.
    \item The indices $\ell$ and $m$ will also be used as functions,
    \beq
    \left\{
        \begin{aligned}
        \ell(\bmu) &= \ell(\mu_1,\mu_2) := \min \{\ell 
        \text{ for which } \bmu \in \cT^\ell \},\\
        m(n) &:= m_n := m 
        \text{ for which } t_n \in [t_{n_m}, t_{n_{m+1}}).
        \end{aligned}
    \right.
    \label{eq:lm}
    \eeq
    For example, given any point $(\bmu,t) \in \cM_T$, we can compute
    $\ell = \ell(\bmu)$ and $m = m_n$ which will allow us to 
    find the parameter-time element the point belongs to, i.e 
    $(\bmu,t) \in \cE^\ell_m.$
    \item Time-partition $[t_{n_m},t_{n_{m+1}})$ for $\cE_m^\ell$ will be 
        given by the following $m_n$ and $n_m$,
    \beq
       \begin{cases}
       m_0 = 0,   \\
       m_n = 1+ \left\lceil \frac{n-1}{20} \right\rceil \text{ if } n > 0,
       \end{cases}
    \quad
       \begin{cases}
       n_0 = 0,   \\
       n_m = 1 + 20(m-1) \text{ if } m > 0.
       \end{cases}
    \eeq
    Apart from the first partition $m=0$, 
    the time partition will contain 20 HFM time-steps.
\end{itemize}
The main motivation for this discretization is to define a basis which is to 
be used by the ROM locally in $\cE^\ell_m$, as we will describe
in \cref{sec:rom}. 

We chose $\cE^\ell_m$ to share the
time-coordinates $\{t_{m_n}\}$ across all parameters $\bmu \in \cM$
for the simplicity of exposition and implementation. But a more flexible set of
elements will be better suited in practice for two important reasons: 
one, to satisfy a stability condition 
(defined in \cref{sec:loc_basis}, \cref{cond:signature}) 
and two, to obtain an optimal number of basis functions through dimensionality
reduction. The construction of elements
that satisfy the two conditions once the parameter space is discretized 
is straightforward: we believe the main 
difficulty lies in the discretization of the parameter space. 
Such difficulties are numerous and challenging, especially in higher
dimensions, but they are not the main focus of this paper. 

\ezfiguretwo{width=0.46\textwidth,trim={0 0 0 50pt}}{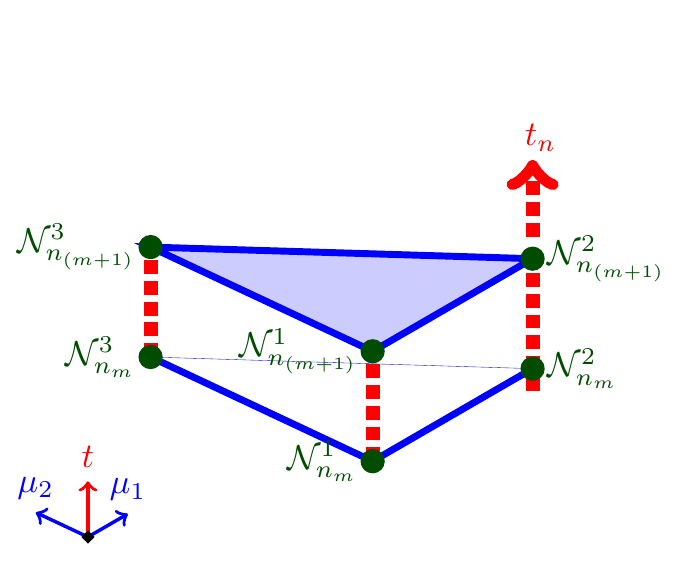}%
            {width=0.46\textwidth,trim={0 0 0 50pt}}{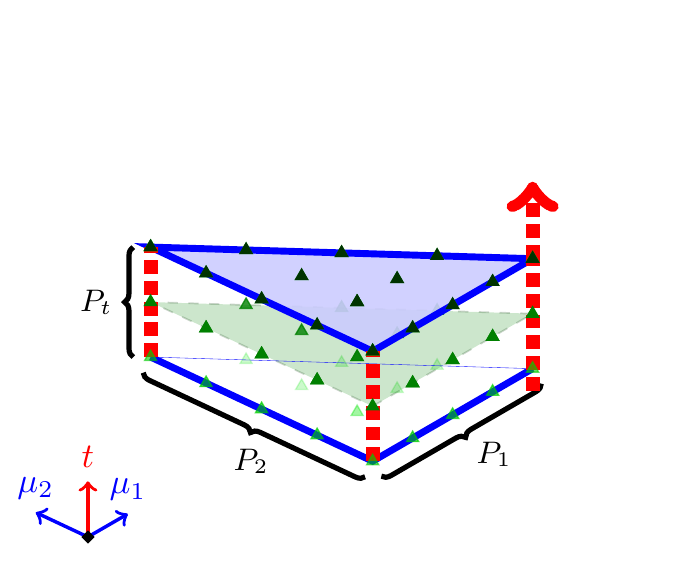}%
            {A diagram of $\cE^\ell_m$ (left) and the 
             uniform interpolation points, as denoted by triangles (right). 
            The dashed lines imply that the time-steps taken by
            the HFM can be finer than the height of $\cE^\ell_m$.}%
            {fig:elem}

\subsection{Displacement interpolation}\label{sec:dinterp}

The notion of \emph{displacement interpolation} is crucial for the success of 
our model reduction. It originally appeared in the optimal transport literature, 
and is also called \emph{McCann's interpolation} since it was devised 
in the study of attracting gas models in the seminal work of McCann 
\cite{mccann97}. 
\emph{Monotone rearrangement} is a simple solution to the Monge-Kantorovich 
optimal transport problem \cite{Kantorovich} between two Borel probability 
distributions over $\RR$.
In our context, the two distributions will be taken simply as two functions 
in $\cC^1(\RR)$. 
The solution is a transport map that deforms one function 
to the other, and displacement interpolation is an interpolation procedure over 
the transport maps, as opposed to interpolation over the functions themselves. 
For a more comprehensive discussion, we refer the reader to the standard 
references \cite{villani03,villani2008}. 

Here, we will follow the definitions and notations of \cite{rim17mr}.
There, the displacement interpolation between two functions of arbitrary sign 
is defined as the monotone rearrangement between positive parts and negative 
parts of the two functions, but we will only need interpolation of 
two non-negative functions (equation (2.18) in \cite{rim17mr}).
\begin{itemize}
    \item The \emph{displacement interpolant} 
    between two non-negative functions $u_1$ and
    $u_2$ by monotone rearrangement will be denoted by,
   \beq
        \cI(u_1,u_2; \alpha) \Ffor 0 \le \alpha \le 1.
    \label{eq:dip_1d}
    \eeq   
   In a setting when $u_1, u_2 > 0$ are both smooth functions 
    in $\cC^\infty(\RR)$,
    \beq
        \cI(u_1,u_2; \alpha) = u_2((1-\alpha)x + \alpha M(x))
        \quad \text{ where } 
        M'(x) = \frac{u_1(x)}{u_2(M(x))},
    \eeq    
    and the map $M(x)$ solves the Monge-Kantorovich problem:
    \beq
        \begin{aligned}
        \text{find } M: \RR \to \RR &\text{ that minimizes } 
        \int_\RR \lvert x - M(x) \rvert u_1(x) \, dx,\\
        \quad
        &\text{ subject to } u_1(M(x)) = u_2(x).
        \end{aligned}
    \eeq
    \item The displacement interpolation between multiple functions is a
    natural extension of \cref{eq:dip_1d}, and will be denoted by,
    \beq
        \cI(u_1,u_2,u_3,u_4; \balpha).
    \label{eq:dip_nd}
    \eeq
    Throughout this paper the parameters $\balpha$ will be three-dimensional,
    \beq
        \balpha = (\mu_1,\mu_2,t),
    \eeq
    and the functions $\{u_j\}$ will be the HFM approximations 
    to \cref{eq:uxtmuN}.
\end{itemize}

Displacement interpolation provides a natural means of obtaining a smooth
deformation between functions whose energy is concentrated in
different spatial locations. When linear interpolation is applied, there is an 
\emph{instant} transfer of energy across large distances, which appears unnatural
for hyperbolic problems. Such behavior is automatically built into standard 
projection-based model reduction, as recognized early on in \cite{rowley00}. 
Coincidentally, precisely this difference in behavior led to
the notion of using interpolation in the transport map rather than
the function itself in optimal transport literature \cite{mccann97}.

Perhaps an example of displacement interpolants that will be familiar to the
reader will be the two-parameter family of displacement interpolants that 
result in  the wavelet basis \cite{daubechies}. Taking a mother-wavelet
along with its dilates and translates as the three functions
and performing displacement interpolation, one obtains the wavelet basis 
functions up to scalar multiples.
    
\section{Construction of local basis}\label{sec:basis}

As the first step in our model reduction, we generate a local basis which 
can represent the solution well locally in the parameter-time space.
Displacement interpolation is an indispensible tool in our development,
and its adaptation will be used to obtain the desired basis functions.

\subsection{Displacement interpolation by pieces and signature}
\label{sec:dinterp_sign}

\ezfigure{width=0.8\textwidth}{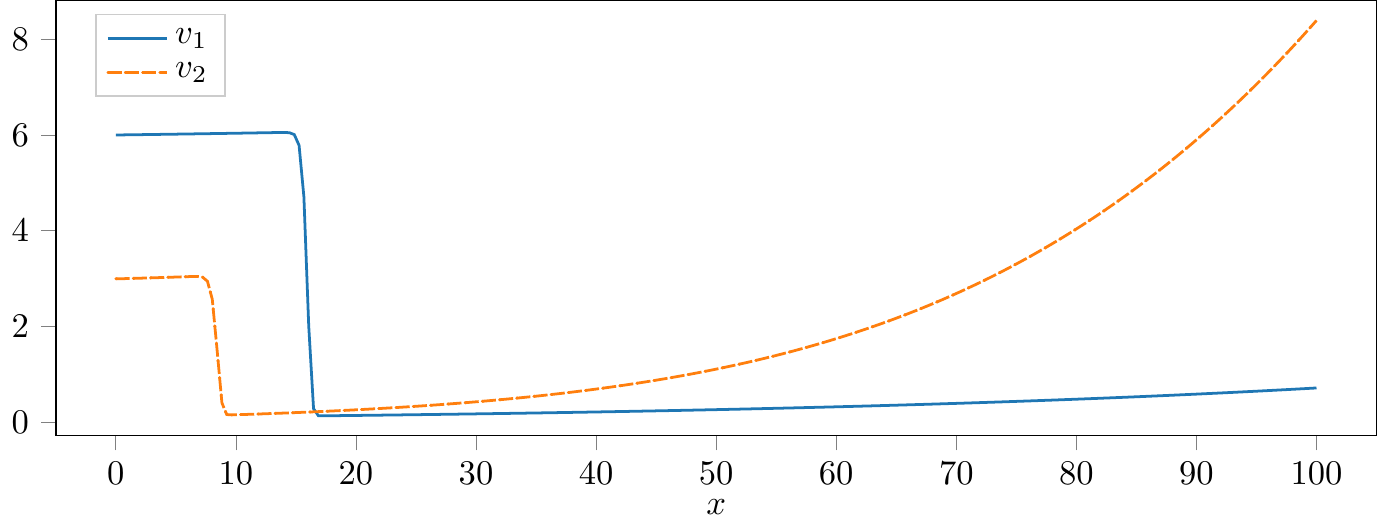}%
         {Two example functions $v_1$ and $v_2$.  
          They are HFM solutions at parameter and time values \cref{eq:uxtmuN},
          $(\mu_1,\mu_2,t) = (3,0.05,5)$ and $(\mu_1,\mu_2,t) = (6,0.02,5)$, 
           respectively.}%
         {fig:two_fctns}%

In this section, we define an adaptation of the usual displacement
interpolation by monotone rearrangement denoted by $\cI$ 
in \cref{eq:dip_1d,eq:dip_nd}.  The functions discussed in this section 
will belong to a function space denoted by $\cU$. We will let $\cU$ be the 
space of piecewise linear functions for simplicity. Even though our 
solutions $(U_{i,n})$ are approximations to solutions to \cref{eq:par_burg} that
may develop jump discontinuities due to the presence of shocks, we will assume 
they can be well represented by piecewise linear functions with sharp
gradients, at least for the purposes of displacement interpolation.
The definitions here extend naturally to Borel measures and distributions, 
but such extensions will not be discussed here.

In our adaptation, we will decompose the derivative $dv /dx$ of a function 
$v \in \cU$, into 
positive and negative pieces, then apply the interpolation $\cI(\cdot)$ between
corresponding pairs of pieces. For example, consider the two
functions in \cref{fig:two_fctns}: the decomposition of its derivative
is shown in \cref{fig:dinterp_parts}.
The motivation for this derivation will be discussed in the following
\cref{sec:loc_basis}.

\ezfigure{width=0.8\textwidth}{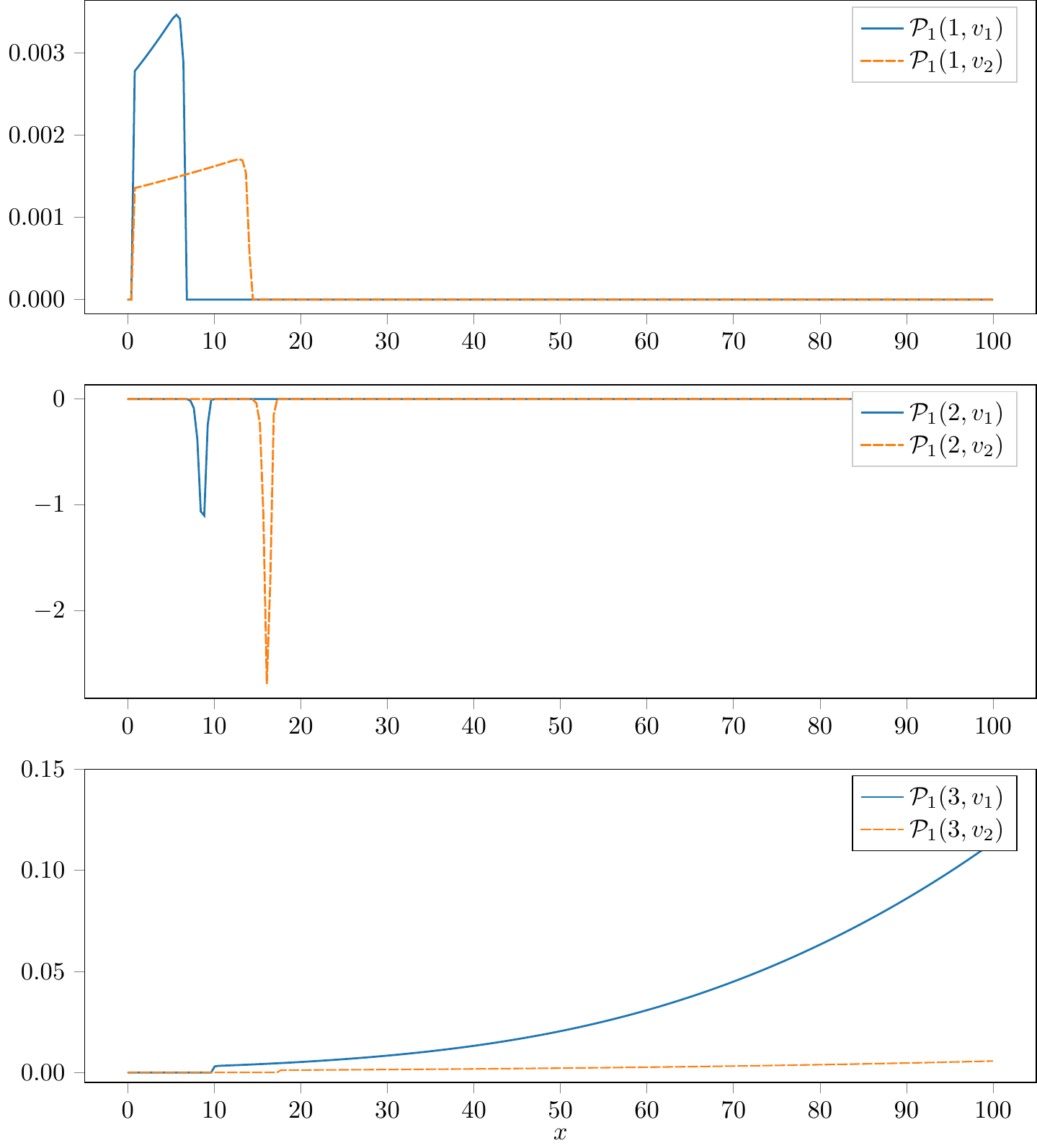}%
         {The pairing of pieces $\cP_2(v_1,v_2) $\cref{eq:P2} for 
          the two functions $v_1$,$v_2$ in \cref{fig:two_fctns} with same 
          signature. Pairs $(\cP_1(1,v_1),\cP_1(1,v_2))$ (top),
          $(\cP_1(2,v_1),\cP_1(2,v_2))$ (middle), and
          $(\cP_1(3,v_1),\cP_1(3,v_2))$ (bottom).
          }%
         {fig:dinterp_parts}%

Now we will explicitly define the procedure outlined above.
\begin{enumerate}
\item Define a set $\cP$, which contains open subsets of the spatial domain,
\beq
\cS(v) := \left\{x : \frac{dv}{dx} (x) \ne 0 \right\},
\quad
\cP(v) := \left\{ \text{connected components of } \cS(v) \right\}.
\label{eq:P}
\eeq
Since $v \in \cC^1(\RR)$, members of $\cS(v)$ are open intervals in which 
the derivative $dv/dx$ has the same sign.
\item Compute the spatial location of these components via the map 
\beq
\mathfrak{c}_v: \cP(v) \to \RR 
\quad \text{ given by } \quad
\mathfrak{c}_v(P) = \frac{1}{\abs{P}} \int_\RR x \cdot \chi_P  \dx,
\eeq
where $\chi_P$ is the characteristic function for the set $P$.
So $\mathfrak{c}_v$ maps each open interval in $\cP(v)$ to its centroid.
\item Denoting the natural ordering of the $\rg(\mathfrak{c}_u)$ by
\beq
\mathfrak{q}_v: \mathfrak{c}(\cP(v)) \to \NN
\label{eq:ordering}
\eeq
we let $\mathfrak{s}_v := \mathfrak{q}_v \circ \mathfrak{c}_v$.
Then $\mathfrak{s}_v$ maps each open interval in $\cP(v)$ to its order 
of position in the domain, as counted from the left. 
\item Next, we define the function which cuts off $dv/dx$ so that it is zero 
outside of the designated interval in $\cP(v)$, 
\beq
\cP_1(\cdot,v) : \rg(\mathfrak{s}_v ) \to \cC(\RR)
\quad \text{ given by } \quad
\cP_1(n,v) := \chi_{\mathfrak{s}_v^{-1} (n)} \, \frac{dv}{dx} .
\eeq
These individual functions will be called \emph{pieces} of $dv/dx$.

\item Now, consider the two functions $v_1$ and $v_2$ that satisfies 
\beq
\abs{\rg(\mathfrak{s}_{v_1})} = \abs{\rg(\mathfrak{s}_{v_2} )}=:J.
\eeq
We define a function that maps $v_1$ and $v_2$ to a set of pairs,
$\cP_2: \cU \times \cU \to 2^{\cU \times \cU}$,
\beq
\cP_2(v_1,v_2) = \bigcup_{j=1}^{J} 
\left\{(w_1,w_2): 
    w_1 = \cP_1(j,v_1) \text{ and } 
    w_2 = \cP_1(j,v_2)
\right\}.
\label{eq:P2}
\eeq
For example, \cref{fig:dinterp_parts} shows an example of the pairs 
in $\cP_2(v_1,v_2)$ for the two functions \cref{fig:two_fctns}.
\item 
Finally, we apply the displacement interpolation $\cI$ to each of the pairs
in $\cP_2(v_1,v_2)$, then compute their cumulative distribution functions (CDFs).
\end{enumerate}

Our strategy is to include these resulting CDFs in our local basis.
The following definition summarizes the above procedure.

\begin{definition}\label{def:dip}
\emph{Displacement interpolation by pieces} for two functions 
$v_1,v_2 \in \cC^1(\RR)$ is given by 
\beq
\cI_P(v_1,v_2; \alpha) 
:= \bigcup_{(w_1,w_2) \in \cP_2(v_1,v_2)} 
\left\{ \int_{-\infty}^x w(z) \dz: w = \cI(w_1,w_2;\alpha) \right\}.
\label{eq:dip}
\eeq
\end{definition}

Multi-dimensional version of $\cI_P$ for $\balpha$ follows straightforwardly 
from the above definition. We will always apply this interpolant on the 
nodes of the element $\cE^\ell_m$ \cref{eq:nodes} so we will make use of
the short-hand notation,
\beq
\cI_P(\cE^\ell_m; \balpha) 
:= \cI_P(u_1,u_2,u_3,u_4; \balpha), 
\quad
\{u_j\} \subset \{u(x,t_*;\bmu_*): (\bmu_*,t_*) \in \cN(\cE^\ell_m) \},
\label{eq:dip_elem}
\eeq
where the subset $\{u_j\}$ is chosen so that $\balpha$ lies in the 
convex hull of the corresponding nodes in $\cN(\cE^\ell_m)$.

The primary motivation for applying the displacement interpolation to 
derivatives \eqref{eq:P} is due to the observation 
that for a class of solutions to hyperbolic conservation laws, 
the spatial derivative $\p u/ \p x$ of the solution evolves along 
a nonlinear vector field that coincides with that constructed by monotone 
rearrangement.
Therefore, it appears reasonable to pursue a construction which interpolates
between derivatives of the solutions. For the full proof of this, 
we refer the reader to \cite{bolley05}.
Since this was proved therein only for non-increasing solutions, we will take 
on a slightly different viewpoint and argue below in \cref{sec:approx} that 
the displacement interpolant is an approximation to an \emph{unknown}
vector field that transforms the solution across the $\{\mu_1,\mu_2,t\}$-axis,
up to first order. Since we are already able to compute the high-fidelity
solution to a desired accuracy, all that we wish to accomplish is to fill 
in the possible states of the solution in between computed solutions 
using the displacement interpolant \cref{eq:dip}.

Another motivation for treating the pieces of the solution with positive and 
negative gradient separately is that for nonlinear conservation laws, the
positive and negative gradients lead to very different dynamics.
In the solutions to \cref{eq:par_burg} during times $0 < t \ll 1$, 
the gradient of the solution $du/dx$ have a positive piece, followed by a 
negative piece and a positive piece. The negative piece corresponds to 
a shock wave, whose propagation is governed by the Rankine-Hugoniot jump 
condition, where the positive pieces correspond to rarefaction 
waves \cite{fvmbook}.
It appears natural to approach the approximation of these behaviors 
separately.

The signs of each piece in \cref{eq:P2} were crucial in determining the pieces 
in \cref{eq:P}. They describe
the qualitative behavior of the function and will serve a role in our stability
condition to be imposed in the next subsection. We will refer to these signs
through the following definition.
\begin{definition}\label{def:sgn}
We define 
$\cS_1:  \cU \, \to \, \{+1,-1\}^{\abs{\rg(\mathfrak{s}_u)}}$,
\beq
\left( \cS_1(v) \right)_n := \sgn(\cP_1(n,v)).
\eeq
We will call $\cS_1(v)$ the \emph{signature} of $v$.
\end{definition}
For example, the signature for the two functions displayed in 
\cref{fig:two_fctns} are,
\beq
\cS_1(v_1) = \cS_1(v_2) = [+,-,+].
\eeq

\subsection{Signature condition and local basis}\label{sec:loc_basis}

In this section, we will construct a basis that is valid for the local element 
$\cE^\ell_m$ using displacement interpolants \cref{eq:dip_elem}.
However, we must first impose some conditions on the local element $\cE^\ell_m$
itself that guarantees that the construction of bases is at all possible.
Consider the functions $\phi_1,\phi_2$ and $\phi_3$ shown in 
\cref{fig:signature}. The three functions each have the signatures
\beq
\cS_1(\phi_1) = [+,-], 
\quad
\cS_1(\phi_2) = [+,-], 
\quad
\cS_1(\phi_3) = [+]. 
\eeq
Since $\phi_3$ has a different signature, it is not an easy task to guess 
the location of the narrow hat function between $\phi_2$ and $\phi_3$,
or between $\phi_1$ and $\phi_3$ for that matter. One approach would be 
to extrapolate the hat function, but this may require significant
effort to be applicable in a general setting.
Instead, we will define a stability condition that prevents such 
difficulties. 
The following condition ensures that the behavior of the solution 
at the nodes of $\cE^\ell_m$ are essentially the same. 

\ezfigure{width=0.8\textwidth}{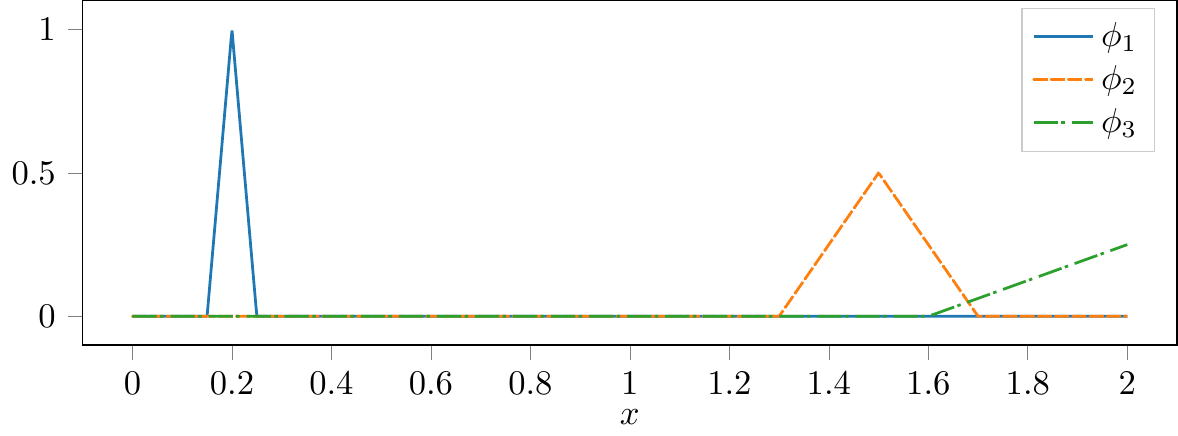}%
         {Three functions $\phi_1,\phi_2,\phi_3$ whose signature $\cS_1$
          is the same for the first two, but different for the last.
          They emulate snapshots of a hat-shaped pulse leaving the domain 
          to the right.}%
         {fig:signature}

\begin{cond}[Signature condition] \label{cond:signature}
Element $\cE^\ell_m$ is said to satisfy the \emph{signature condition} if
each $\balpha \in \cE^\ell_m $ lies in the convex hull of some 
$\cN_0 \subset \cN(\cE^\ell_m)$ such that $\cS_1(\cN_j) = \cS_1(\cN_k)$
for all $\cN_j,\cN_k \in \cN_0$.
\end{cond}

\cref{cond:signature} can be also seen as a restriction
on how much you can expand each element $\cE^\ell_m$ without running into a
significant change in the solution. 

Now, one natural consequence of this condition is that a more flexible 
form of $\cE^\ell_m$ will be required to fully take advantage of the 
displacement interpolation. 
Consider the behavior of our solutions displayed in \cref{fig:hfm_sols}.
The solution is initially $u(x,0) = 0$ for all parameters, so the 
solution has null signature at $t =0$. Then there is an immediate change of 
signature for the solution at positive time to $[+,-,+]$, due to the incoming 
boundary condition to the left and the exponential source term 
\cref{eq:par_burg}.
So $\cE^\ell_0$ and $\cE^\ell_1$ must separate the initial condition and the
first time-step. Depending on the parameter, the shock travels to the right and 
leaves the domain at different times depending on the parameters. Therefore
to maximize the parameter-time regions covered by $\cE^\ell_m$, one must allow
an unstructured discretization. This has the implication that 
the final time-step that could be allowed for the reduced-order model would 
be constrained by adhering to the uniform-in-time discretization in
\cref{eq:cE}. Nonetheless, we will maintain the structured form of $\cE^\ell_m$ 
for the simplicity of exposition and implementation.

As parameters and time are varied in the interpolation procedure
$\cI_P$ defined above \cref{eq:dip},
the location and the sharpness of the negative gradient representing the shock 
are also varied smoothly. Furthermore in \cref{eq:dip} we are 
constructing multiple interpolants corresponding to each pair in 
$\cP_2(v_1,v_2)$ \cref{eq:P2} separately, then including it into the set of 
functions on the RHS in \cref{eq:dip}.
This allows more flexibility in the basis representation: for an unknown
parameter value, the shock may be in a slightly different location than the
computed interpolant, therefore we want to ensure that the local basis 
is able to represent local perturbations of the interpolant as well. 
This may increase the number of basis functions,
but the dimensionality reduction is secondary at this stage:  
the goal of our basis generation is to include 
all basis functions that may potentially be used to represent the solution.
 Once we have generated such bases, we will be
able to further extract a reduced basis with much smaller dimensionality
(see \cref{sec:pod_proj} below).

Now we are ready to construct the local basis.
We sample $\cE^\ell_m$ uniformly at points $\{\balpha_k \in
\cE^\ell_m\}$ that are equally spaced along each $\{\mu_1$,$\mu_2,t\}$-axis 
(i.e, a restriction of the Cartesian grid onto $\cE^\ell_m$).
We will use $P_t$,$P_1$,and $P_2$ sample points along the $t$-axis, 
$\mu_1$-axis and $\mu_2$-axis, respectively, as illustrated in \cref{fig:elem}.
Using these sets of functions, we construct a local basis 
$\cW^\ell_m$ corresponding to the element $\cE^\ell_m$ by producing 
an orthonormal basis for the linear 
space spanned by the interpolants  $\cI(\cE^\ell_m; \balpha_k)$ using, 
for example, the Gram-Schmidt process,
\beq
\left\{
\begin{aligned}
\Span \cW^\ell_m &= \Span\{\cI(\cE^\ell_m; \balpha_k)\},  \\
\langle w, w \rangle &= 1,
\quad \text{ for } w \in \cW^\ell_m ,\\
\langle w_1, w_2 \rangle &= 0,
\quad \text{ for } w_1,w_2 \in \cW^\ell_m \text{ and } w_1 \ne w_2.
\end{aligned}
\right.
\label{eq:ortho_basis}
\eeq
where $\langle \cdot, \cdot \rangle$ denotes the usual $L^2$ inner product.
The local bases $\cW^\ell_m$ will be used for Galerkin projection in 
\cref{sec:rom}.

We end this section with the remark that the local basis constructions can be
performed in an embarrassingly parallel manner.

\section{Approximation properties of the displacement interpolant}
\label{sec:approx}

As mentioned above, displacement interpolants are often more useful than 
linear interpolants for approximating wave-like phenomena. This is not so 
surprising upon closer inspection, as there are close
relations between fluid equations such as compressible Euler's equation and
the Monge-Kantorovich problem \cite{villani03}.
In this section, we briefly discuss the approximation properties of 
the displacement interpolation $\cI$ \cref{eq:dip_1d} and its
variant $\cI_P$ introduced in \cref{eq:dip}. 

Upon discretization, $\cI$ can be seen as a type of a Lagrangian method 
such as the particle-in-cell (PIC) or large-time-step (LTS) 
methods. 
However, unlike in most applications of these methods, the
velocity of the particles depend on the \emph{dependent variable} 
rather than the spatial variable $x$: particles $x(t)$ can travel at different 
speeds although they are passing through the same spatial location, as evident 
in \cref{fig:vector_field}. In other
words, the displacement interpolant represents the dynamics of a wave profile 
moving through a \emph{time-dependent} or \emph{nonlinear} medium \cref{eq:pde}.
In fact, we will show that the interpolant $\cI$
is a first-order approximation of a smooth time-dependent and possibly nonlinear 
velocity field in terms of the spatial variable $x$. On the other hand, it is 
also a \emph{linear} interpolant when viewed as a function of the 
dependent variable (see \cref{eq:x_sum,eq:x_approx} below). 

\ezfigure{width=0.8\textwidth}{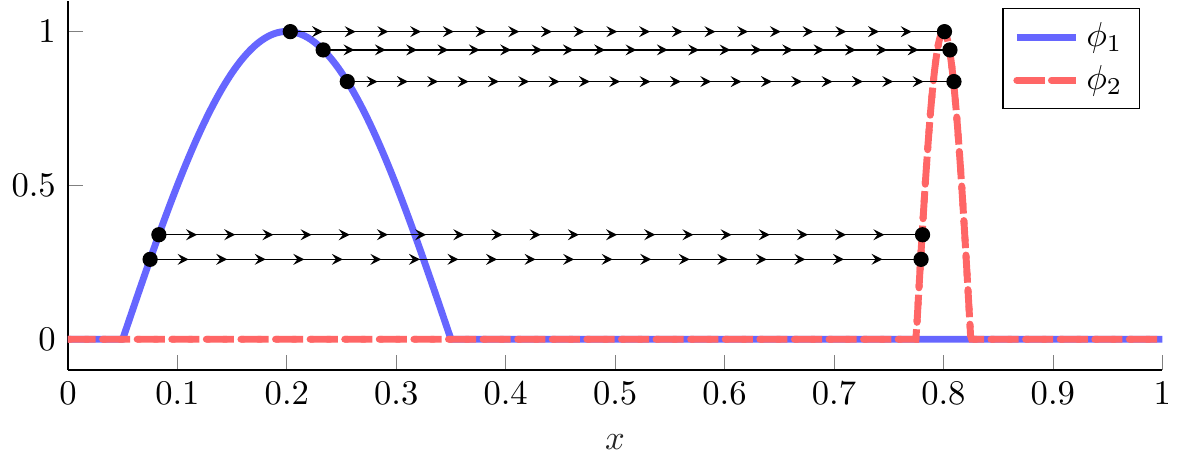}%
         {Particle trajectories for the displacement interpolant
          \cref{eq:dip}.}%
         {fig:vector_field}

To elaborate, let us write both linear interpolation and displacement
interpolation $\cI$ as a solution to a time-dependent problem. 
Here time is seen merely as a parameter, and is not 
necessary the time variable of the underlying differential equation.
Suppose that $u$ is the solution to the problem,
\beq
\frac{\partial u}{\partial t}
+ \frac{\partial}{\partial x} \left( c(x,t) u \right) = \psi(x,t),
\eeq
and further suppose that we know the solution $u$ at two different times
$t_1$ and $t_2$, 
\beq
u(x,t_1) = u_1(x) 
\quad \text{ and } \quad
u(x,t_2) = u_2(x).  
\eeq
We will show that one obtains the linear interpolant if one assumes $c\equiv 0$, 
and the displacement interpolant if one assumes $\psi \equiv 0$.
Throughout, we will assume that $c$ and $\psi$ are smooth functions of time.

\subsection{Linear interpolation} 
Let us first assume $c \equiv 0$, then $u$ would satisfy,
\beq
\frac{\partial u}{\partial t} = \psi(x,t),
\quad \text{ so that } \quad
u(x,t) = u_1(x) + \int_{t_1}^{t} \psi(x,s) \ds.
\eeq
If $\psi(x,\cdot) \in \cC^1([t_1,t_2])$ we can write a first order 
approximation of the integral above,
\beq
\int_{t_1}^{t} \psi(x,s) \ds
=  \frac{u_2(x) - u_1(x)}{t_2 - t_1} (t - t_1) + \cO(\Delta t),
\quad \text{ where } \abs{t_2 - t_1} \le \Delta t,
\label{eq:uint_approx}
\eeq
and we arrive at the linear interpolant
\beq
\bar{u}(x,t) = a_1(t) u_1(x) + a_2(t) u_2(x),
\quad
\text{ where }
a_1(t) = \frac{t_2 - t}{t_2 - t_1}, \quad a_2(t) = \frac{t - t_1}{t_2 - t_1}.
\label{eq:u_approx}
\eeq
Since $\abs{u(x,t) - \bar{u}(x,t)} \le \cO(\Delta t)$,
$\bar{u}$ is a first order approximation to $u$ in $t$-variable.

\subsection{Displacement interpolation} 
Let us now assume $\psi \equiv 0$, so that $u$ solves the conservation law,
\begin{equation}
\frac{\partial u}{\partial t}
+ \frac{\partial}{\partial x} \left( c(x,t) u \right) = 0.
\label{eq:pde}
\end{equation}

Let us assume $u > 0$, then define the CDF of $u(x,t)$,
\[
U(x,t) = \int_{-\infty}^x u(z,t) \dz
\quad
t_1 \le t \le t_2.
\]
Then we can rewrite the PDE in terms via the transformation
$(x,t)  \leftrightarrow (y,\tau)$ where  $y=U(x,t),  \tau = t.$
Under this transformation the derivatives in \eqref{eq:pde} now becomes
\begin{equation}
\begin{aligned}
\frac{\partial}{\partial x} &= 
\frac{\partial U}{\partial x} \frac{\partial}{\partial y} 
= u(x,t) \frac{\partial}{\partial y}, \\
\frac{\partial}{\partial t} &= \frac{\partial}{\partial \tau}
+ \frac{\partial U}{\partial t} \frac{\partial}{\partial y}
= \frac{\partial}{\partial \tau}
-c(x,t) u(x,t) \frac{\partial}{\partial y},
\end{aligned}
\end{equation}
where the last identity comes from,
\begin{equation}
\begin{aligned}
\frac{\partial U(x,t)}{\partial t}
= \frac{\partial}{\partial t} \int_{-\infty}^x u(z,t) \dz
= \int_{-\infty}^x \frac{\partial}{\partial t} u(z,t) \dz 
= - c(x,t) u(x,t).
\end{aligned}
\end{equation}
Expressing the above transformation in variables $y$ and $\tau$,
\beq
v(y,\tau) := \frac{1}{u(x(y,\tau),\tau)} 
\quad \text{ and } \quad
c(y,\tau) := c(U^{-1}(y,\tau),\tau),
\eeq
the transformation can be summarized as,
\begin{equation}
\frac{\partial}{\partial x} = \frac{1}{v(y,\tau)} \frac{\partial}{\partial y},
\qquad
\frac{\partial}{\partial t} = 
\frac{\partial}{\partial \tau} -
\frac{c(y,\tau)}{v(y,\tau)} \frac{\partial}{\partial y}.
\end{equation}
The transformed PDE is
\begin{equation}
\left( \frac{\partial}{\partial \tau} - \frac{c}{v} \frac{\partial}{\partial y}
\right) \left( \frac{1}{v} \right)
+ \frac{1}{v} \frac{\partial}{\partial y} \left( \frac{c}{v} \right) = 0.
\end{equation}
Using the chain rule, and multiplying both sides by $v^2$,
\begin{equation}
\frac{\partial v}{\partial \tau} = \frac{\partial c(y,\tau)}{\partial y}.
\end{equation}
Therefore, the evolution of $v(y,\tau)$ now solely involves the integration of
$c(y,\tau)$. Integrating with respect to $y$ and assume $c \to 0$ as 
$y \to 0$ 
\begin{equation}
\frac{\partial x}{\partial \tau} = c(y,\tau) 
\quad \text{ so that } \quad
x(y,\tau) = x(y,\tau_1) + \int^\tau_{\tau_1} c(y,s) \ds .
\label{eq:dyn}
\end{equation}
where, $\tau_1 = t_1$ and $\tau_2 = t_2$.
Note that this is simply a variant of (9) in \cite{bolley05}, if we assume that 
$c(y,t)$ is in fact just a function of $y$. 

This is precisely the evolution equation for the particles, whose change 
in speed due to the field $c(x,t)$ at $x$ is now implicit in the moving 
coordinate $y = U(x,t)$. As we have observed above, displacement
interpolation shifts the particles at a constant speed, and
each speed is a function of $y$, as opposed to the original spatial variable
$x$. Therefore, the displacement interpolant \eqref{eq:dip_1d} is simply an 
approximation to the dynamics given by \eqref{eq:dyn}, in which the integral is 
approximated by a term linear in $t$.

Let $x_1(y) = x(y,\tau_1)$ and $x_2(y) = x(y,\tau_2)$,
and if $c(y,\cdot) \in \cC^1([\tau_1,\tau_2])$ we can write a first order 
approximation just as in \cref{eq:uint_approx},
\beq
\int_{\tau_1}^{\tau} c(y,s) \ds
=  \frac{x_2(y) - x_1(y)}{\tau_2 - \tau_1} (\tau - \tau_1) + \cO(\Delta \tau),
\quad \text{ where } \abs{\tau_2 - \tau_1} \le \Delta \tau,
\label{eq:xint_approx}
\eeq
The first order approximation is then simply,
\beq
\bar{x}(y,\tau) = a_1(\tau) x_1(y) + a_2(\tau) x_2(y),
\quad
\text{ where }
a_1(\tau) = \frac{\tau_2 - \tau}{\tau_2 - \tau_1}, 
a_2(\tau) = \frac{\tau - \tau_1}{\tau_2 - \tau_1}.
\label{eq:x_approx}
\eeq

Therefore, when $\cI$ is applied to the snapshots from the PDE \eqref{eq:pde}, 
it yields an approximation to the dynamical system \eqref{eq:dyn},
using a first order approximation to the \emph{speed} of the evolution
\cref{eq:xint_approx}. This results in a particle-type method in 
which each particle moves at constant speed over time.

\subsection{Generalized form} \label{sec:gen_form}
There is a clear analogy between \cref{eq:u_approx} and \cref{eq:x_approx}
and the two approximations leading to a more general form.
For the former, we have the usual linear approximation by separation
of variables,
\beq
u(x,t) \approx \bar{u}(x,t) = \sum_{n=1}^\infty a_n(t) f_n(x),
\label{eq:u_sum}
\eeq
for some suitable basis functions $\{f_n\}$, and
for the latter, we have the approximation of the transport map that transforms
some fixed function $u_0(x)$,
\beq
u(x,t) = u_0(x(y,t)) 
\quad \text{ where } \quad
x(y,t) \approx \bar{x}(y,t) = \sum_{n=1}^\infty b_n(t) g_n(y),
\label{eq:x_sum}
\eeq
again for some basis $\{g_n\}$. For example, any linear combination of  
interpolation by pieces $\cI_P$  \cref{eq:dip} could be rewritten in the form
of \cref{eq:x_sum}. 
In effect, we have a separation of variables for the 
transport map. One may also view it as a series representation of the 
characteristics equations. 

The evolution equation \cref{eq:dyn} is closely related to the Wasserstein
metric in optimal transport \cite{bolley05}, however the connection is no
longer obvious for the more general form of the approximation above in 
\cref{eq:x_sum}. In other words, solution to the Monge-Kantorovich problem 
appears naturally as the first order approximation for the equation 
\cref{eq:pde}, but the higher order approximations are no longer solutions to
this problem, as they do not move each unit of mass with constant speed,
nor is it transport along multiple rays \cite{caffarelli02}.

This viewpoint is already visible in the previous perspectives presented in
in \cite{welper17,welper17p,rim17reversal,schulze18}. A common thread among 
them is that they naturally lead to an optimization problem. The general form 
here seems to offer a more constructive insight, and relates the resulting 
interpolation to classical approximation theory, simply by posing the 
problem as an approximation problem in the $y$-variable.
This general form also can be extended to multi-dimensional settings where
$x \in \RR^n$ via  the Radon transform \cite{radonsplit,Helgason2011}, since 
the approximation above can be applied to each direction $\omega \in S^{n-1}$
as discussed in \cite{rim17mr}.
In any case, the formulation \cref{eq:x_sum} perhaps will serve to complement 
standard tools in linear approximation, such as the SVD.

\section{Galerkin projection}\label{sec:rom}
Once the local basis has been constructed as described in \cref{sec:basis}, 
we are ready to set up the Galerkin projection. The projection itself does
not require special treatment, other than the fact that the basis functions
will change locally with respect to parameters $(\mu_1,\mu_2)$ and time $t$. 
Throughout this section, we will make use of
Einstein's summation notation, in which repeated indices are summed over.

A given pair of parameter values $(\mu_1,\mu_2)$ will 
belong to an element $\cT^\ell$ \cref{eq:tri} 
and the $n$-th time step $t_n$ to an interval $[t_{m_n}, t_{m_n+1})$,
so we will employ the notations \cref{eq:lm},
$\ello = \ell(\mu_1,\mu_2)$ and $m = m_n.$
The local basis will be denoted by $\cW^\ell_m$ \cref{eq:ortho_basis},
and will be valid locally in the parameter-and-time element 
$\cE^\ell_m$ \cref{eq:cE}. In our numerical discretization, the basis 
will be represented by a matrix $\bfW^\ell_m$ of size $\RR^{N\times M^\ello_m}$
\beq
\left(\bfW^\ell_m \right)_{ij}= 
 W_{ij,m}^\ell 
\quad \text{ where } \quad 
1 \le i \le N, 
\quad
1 \le j \le M^\ell_m.
\label{eq:matrix_basis}
\eeq
The matrix will be orthogonal, that is, 
$W_{ij,m}^\ell W_{jk,m}^\ell = \delta_{ik}.$
We will be representing the finite volume solution 
$U^\ell_{i,n}$ \cref{eq:fvcell} in this local basis by $r_{j,n}^\ell$,
\beq
r^\ello_{j,n} = W_{ij,m_n}^\ello U_{i,n}^\ello,
\quad \text{ and } \quad
U^{\ello}_{i,n} = W^{\ello}_{ij,m_n} r^{\ello}_{j,n}. 
\eeq
One must also compute the \emph{transition matrix} whenever time-stepping
between $\cE^\ell_m$ and $\cE^\ell_{m+1}$. That is, when $m_{n+1} \ne m_n$,
\beq
T_{ij,m_{n+1}} = W_{ik,m_{n+1}}^\ell W_{kj,m_n}^\ell,
\quad \text{ and } \quad
\tilde{r}_{i,n}^\ell = T_{ij,m_{n+1}} r_{j,n}^\ell
\eeq
which should then be used for the initial calculations in $\cE^\ell_{m+1}$, 
but the details are straightforward and we will omit the specifics here.

\subsection{Reduced order model} 

In this subsection, we will restrict our attention to
a fixed parameter-time element $\cE^\ell_m$ and suppress the indices
$\ell=\ell(\mu_1,\mu_2)$ and $m= m_n$ for the ease of notation by letting
\beq
W_{ij} = W_{ij,m_n}^\ell,
\quad
U_{i,n} = U_{i,n}^\ell
\quad \text{ and } \quad
r_{j,n} = r_{j,n}^\ell.
\eeq

We project the finite volume update \cref{eq:hfm_update} to the local basis
$(\bfW)_{ij} = W_{ij}$ of the size $\bfW \in \RR^{N \times M}$. 
This yields the time-update in the new basis,
\beq
    \begin{aligned}
    W_{ij} U_{i,n+1} &= W_{ij} U_{i,n}  \\
        & \qquad - \frac{\Delta t}{\Delta x}(W_{ij}
                   [U_{i+1,n}]^2 - W_{ij}[U_{i,n}]^2)
                + k \, W_{ij} \left( 0.02 e^{\mu_2 x_i} \right).
    \end{aligned}
    \label{eq:rom_update}
\eeq
The last two terms on the right still require computations
that depend on the degree of freedom of the HFM, so we must treat them
separately for performance.
\begin{itemize}
\item To compute the projection of the flux on to the new basis, we 
      simply need 
\beq
\begin{aligned}
        W_{ij} [U_{i,n}]^2 
            &= W_{ij} (W_{ki}r_{k,n}) 
                    (W_{p i}r_{p,n}) \\
          &= \underbrace{(W_{ij} W_{ki} W_{pi})}%
              _{M \times M \times M} r_{k,n} r_{p,n}
=: F_{jkp} r_{k,n} r_{p,n}.
\end{aligned}
\label{eq:rom_flux}
\eeq
So the tensor $F_{jkp}$ can be computed during the offline
stage. This is due to the fact that the flux function $f(u) = u^2/2$ is a 
low-order polynomial. Fortunately, many common nonlinear hyperbolic problems 
have flux functions that are low-order polynomials, allowing us
to employ this simple and exact reduced-basis representations of the nonlinear
terms. 
\item To compute the projection of the source term onto these basis functions, 
we make use of a truncated Taylor series expansion of the exponential function
up to $Q$-th term. Then we arrive at the approximation,
\beq
W_{ij} \left( 0.02 e^{\mu_2 x_i} \right)
\approx 
0.02 \underbrace{W_{ij} \frac{1}{q!}(x_{i})^q}_{M \times Q}
(\mu_2)^q  =: S_{jq} (\mu_2)^q.
\label{eq:rom_src}
\eeq
The $S_{jq}$ in the RHS can be precomputed offline, and 
updated whenever $\mu_2$ is determined. During the online-stage 
one can compute
\beq
(1,\mu_2, \mu_2^2 , \cdots, \mu_2^q , \cdots , \mu_2^Q),
\eeq
then the dot product in \cref{eq:rom_src} can be computed with small cost.
\end{itemize}
In the case the flux terms or source terms are not in polynomial 
form \cref{eq:rom_flux} or are not easy to approximate as functions of the 
parameters \cref{eq:rom_src}, one must make use of more 
advanced numerical techniques, such as the discrete empirical interpolation 
method (DEIM) \cite{chanturantabut10}, 
Gauss–Newton with approximated tensors (GNAT) \cite{carlberg13},
or local DEIM \cite{peherstorfer14}, for example.

Finally, we obtain the time-update in the local basis,
\beq
 r_{j,n+1} = r_{j,n} 
                   - \frac{\Delta t}{\Delta x} F_{jkp} r_{k,n} r_{p,n}
                   + \Delta t S_{jq} (\mu_2)^q.
\label{eq:rrom_update}
\eeq
The system dimension for each update depends only on $M,Q \ll N$. 

The computational cost for the time-stepping in the new basis now scales 
linearly, requiring $\cO(N)$ floating-operations per
one solve. In higher spatial dimensions the performance gain is expected 
to be more substantial. Nonetheless, the number of time-steps taken  
still depends on the dimension of the HFM and this burden is now the bottleneck 
of the ROM, although there may well be further model reduction possible
for the time variable: this will be a topic of future investigation.

In this paper, we will compute the ROM solution
\cref{eq:rrom_update} up to time $t = 12$ for ease of implementation
in order to satisfy \cref{cond:signature} with our structured discretization
of the parameter-time space \cref{eq:cE},
but this is not a restriction on the method itself, and can be lifted through 
an unstructured discretization of $\cE^\ell_m$ as discussed above 
in \cref{sec:loc_basis}.

\subsection{POD projection} \label{sec:pod_proj}
It is possible to reduce dimensionality of the system \cref{eq:rrom_update} 
significantly further. Note that we
have not so far made use of standard dimensionality reduction techniques
such as the proper orthogonal decomposition (POD). 
We do so here by running the ROM \cref{eq:rrom_update} over many
parameter values in $\cM$. Then one can collect snapshots of the solution
in the basis $\cW^\ell_m$, and apply the standard dimensionality reduction
tools to obtain a reduced basis. This reduction does not require additional
computations from the HFM, and can be done completely offline.
Moreover, the elements $\cE^\ell_m$ can be refined further during the 
computations to optimize the number of bases per time-step.  
The rule of thumb is that for hyperbolic problems, the number of basis required 
for accuracy increases linearly with the Euclidean distance in the 
parameter-time space, therefore more refinements will lead to a futher reduction 
in the number of basis functions required per time-step in \cref{eq:rrom_update}.

The number of bases we have obtained with unrefined-and-structured elements
$\cE^\ell_m$ \cref{eq:cE} with the standard POD using the truncation
threshold for the ratio of singular values $\sigma_n / \sigma_1$ 
set to \texttt{1e-8}, are plotted in \cref{fig:basis_no}. 
We also remark that the number of bases for each element $\cE^\ell_m$
appears to be independent of the size of the HFM, as long as the number of 
HFM time-steps included in each $\cE^\ell_m$ \cref{eq:cE} remains constant.
More elements will be required to reach the final time, but the number of local 
bases functions needed at each time-step remains small, maintaining the linear 
scaling $\cO(N)$.

The linear growth of number of bases observed in the number of basis 
for $\cE^\ell_m$ in \cref{fig:basis_no} agrees with the rule
of thumb, and can be controlled at will by the refinement of $\cE^\ell_m$.
For example, if we restrict our attention to 
a sub-region of the parameter-time space, given by
\beq
  \tilde{\cE}_m := \{(\mu_1,\mu_2,t): \mu_1 \in [6,7] , \mu_2 \in [0.06, 0.075], 
                             t \in [t_{n_m},t_{n_{m+1}})\}
        \subset \cE_m^7,
\eeq
the number of reduced bases required to maintain the same of level of
accuracy is dramatically reduced to less than $13$, as plotted 
in \cref{fig:basis_no}. We will continue to employ the originally defined 
elements $\cE^\ell_m$ here, solely for ease of illustration.

\ezfiguretworow{width=0.95\textwidth}{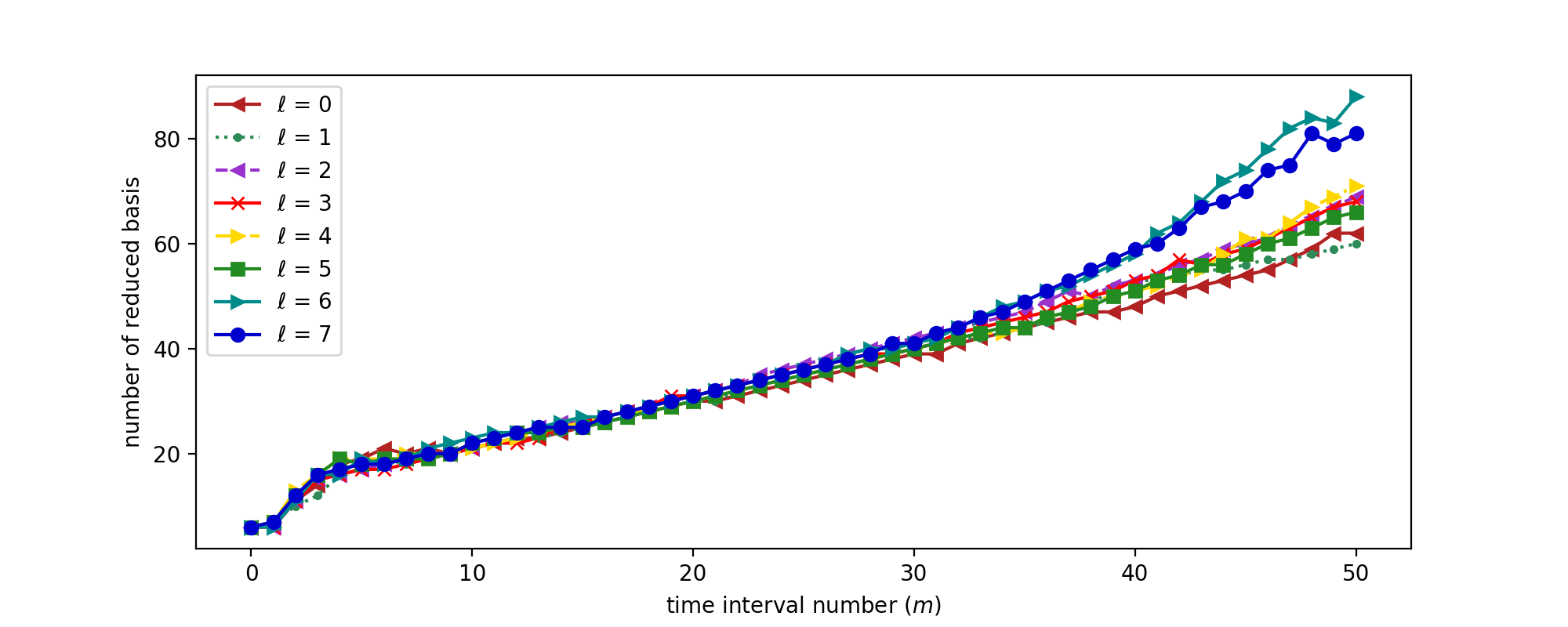}
               {width=0.95\textwidth}{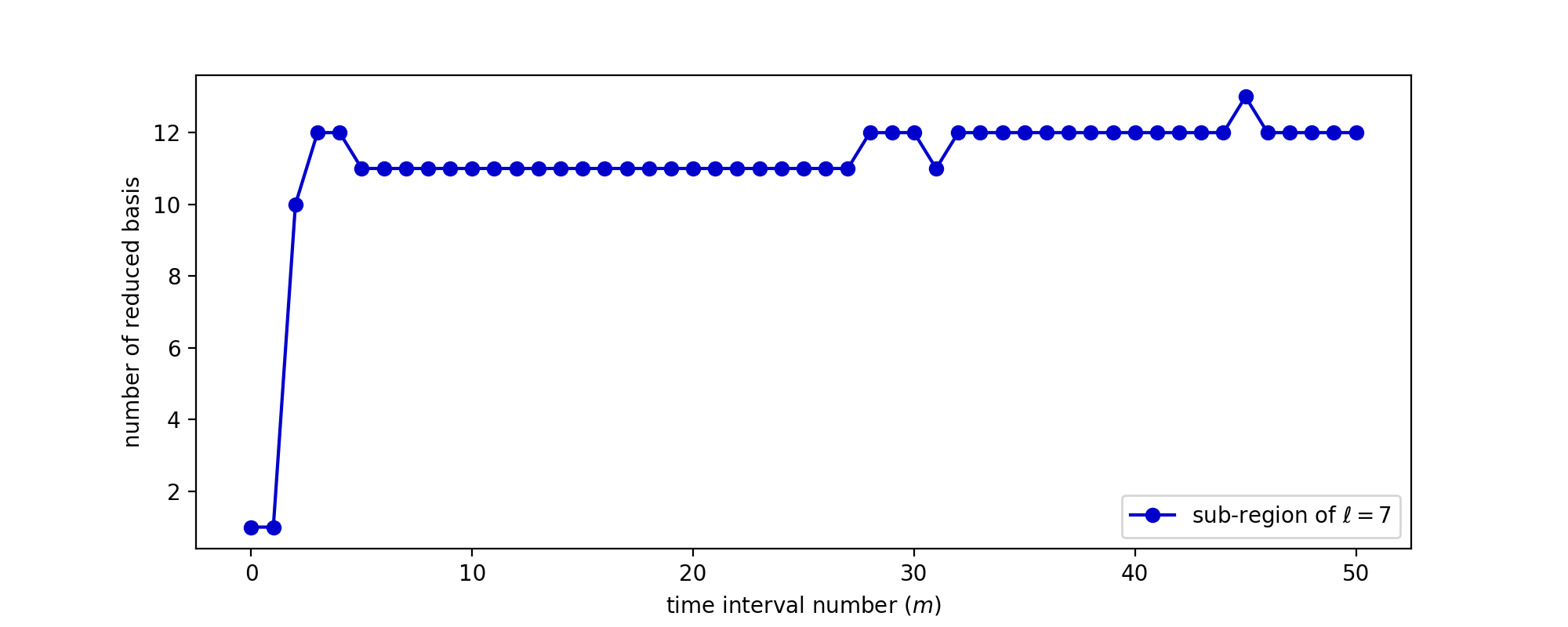}
         {Number of reduced basis for each element $\cE^\ell_m$ (top) and
          number of further reduced basis for in the refined sub-region
          $\tilde{\cE}_m$ lying in the element $\cE^7_m$ (bottom).}%
         {fig:basis_no}


\subsection{Error} 
A comparison of a ROM solution for a randomly chosen parameter value
with its corresponding HFM solution is 
shown in \cref{fig:rom_err} up to final time $t = 12$. 
The HFM differs from the ROM with maximum pointwise relative error of roughly 
\texttt{1e-3}. The error is localized near the shock, which is to be expected.

To estimate the global error, we define the maximum point-wise
relative error as,
\beq
E_{\textrm{Rel}}(\bmu) = 
\max_{t \in (0,12]}
\max_{x \in [0,100]}
        \frac{|u_\textrm{HFM}(x,t;\mu_1,\mu_2) 
             - u_\textrm{ROM}(x,t;\mu_1,\mu_2)|}%
             {|u_\textrm{HFM}(x,t;\mu_1,\mu_2)|}.
\label{eq:rele}
\eeq
Assuming a uniform distribution over the parameter space $\cM$,
we estimate the mean and the variance of the random variable 
$E_{\textrm{Rel}}(\bmu)$. We use Monte Carlo sampling to obtain $10,000$
samples in $\cM$ and compute solutions from both the ROM and the HFM,
then use them to obtain the relative error \eqref{eq:rele}. 
The mean and the variance 
are estimated as follows, indicating good accuracy:
\[
 \mathbb{E}[E_\textrm{Rel}(\bmu)] = \texttt{3.1402e-03}, 
 \quad \text{ and } \quad
 \textrm{Var}[E_\textrm{Rel}(\bmu)] =  \texttt{1.7852e-06}.
\]

\ezfigure{width=0.9\textwidth}{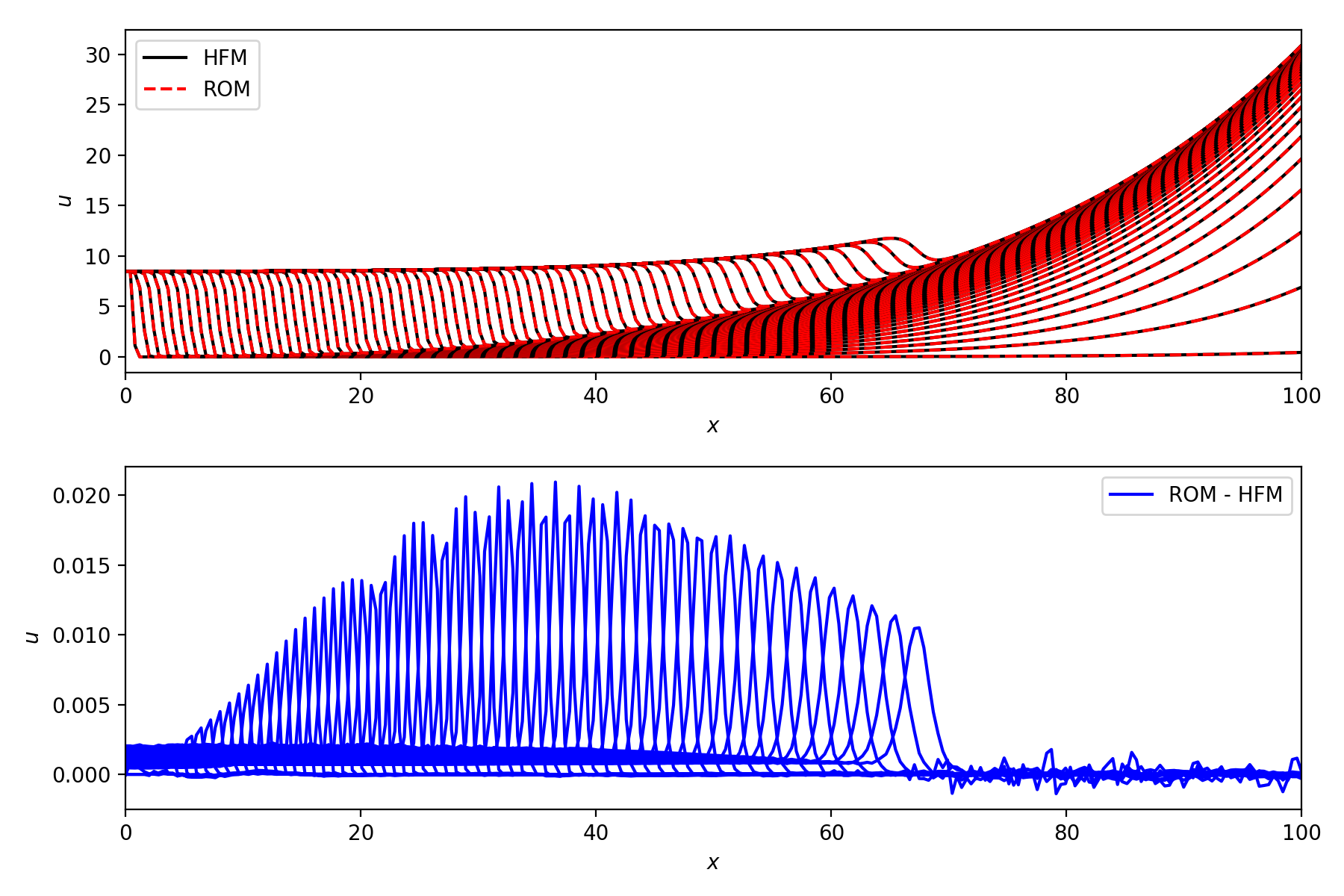}%
            {Comparison between ROM and HFM (top) and the difference
             between the two (bottom) for parameter values
            $(\mu_1,\mu_2) = (8.4601,0.0750)$.
            Both are plotted between every $15$ time-steps. }%
            {fig:rom_err}

\section{Uncertainty quantification}\label{sec:uq}

Using the ROM developed in the previous sections, we perform a few UQ tasks in
this section. We sample the parameter space uniformly, by drawing samples from
\beq
\bmu \sim \cU(\cM).
\eeq
Using standard Monte-Carlo samples, we can easily run the ROM and extract
the quantity of interest (QoI). We sampled 10,000 values of $\bmu$, then ran
the ROM for each value. 

\begin{itemize}
\item We compute the average solution and its variance, for each triangle 
      $\cT^\ell$ in the parameter space. That is, we estimate
      \beq
        \mathbb{E}_{\bmu}[u(x,t;\bmu)]
        \quad \text{ and } \quad
        \mathrm{Var}_{\bmu}[u(x,t;\bmu)].
      \eeq  
      These are plotted in 
      \cref{fig:u_mean_var,fig:u_mean_var_time}.
      The large variance is near the shock since the shock location
      itself is sensitive to the change in parameter, and the large gradient
      near the shock incurs large variance even when the shock location
      varies subtly.

\ezfiguretworow{width=0.9\textwidth}{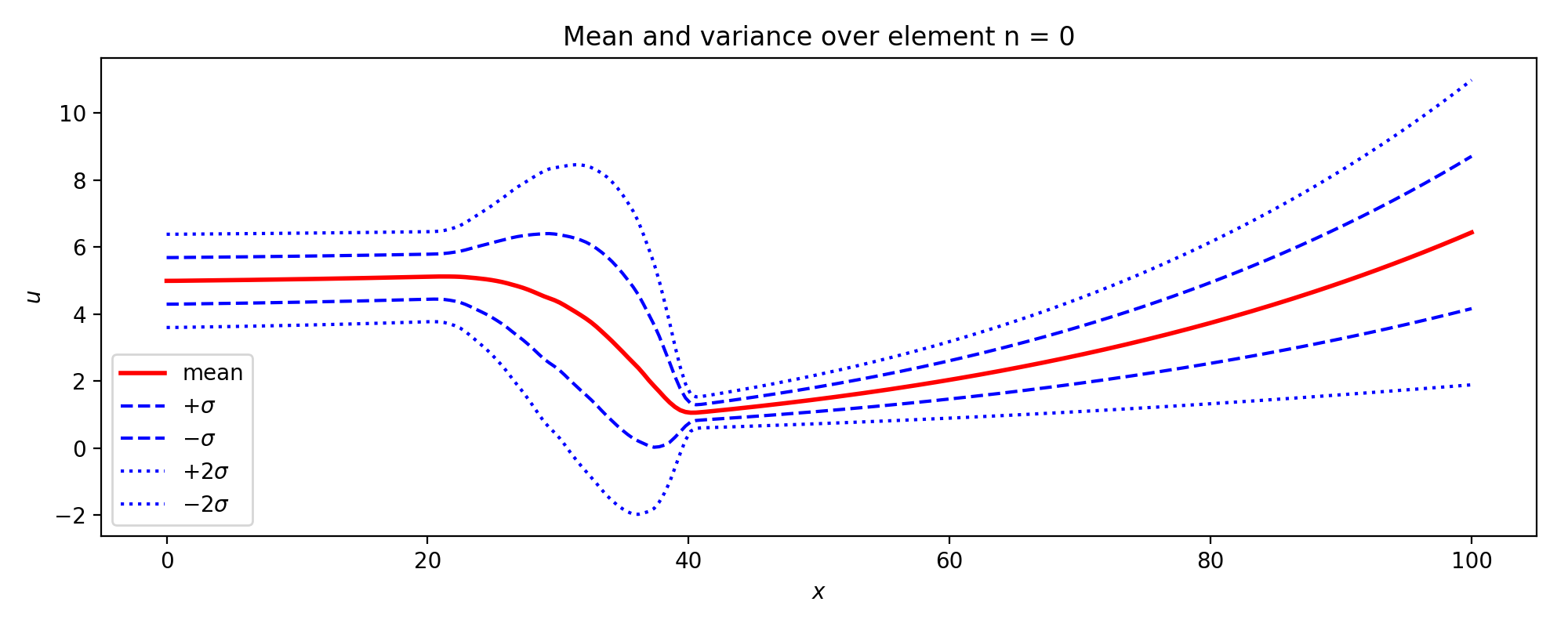}%
               {width=0.9\textwidth}{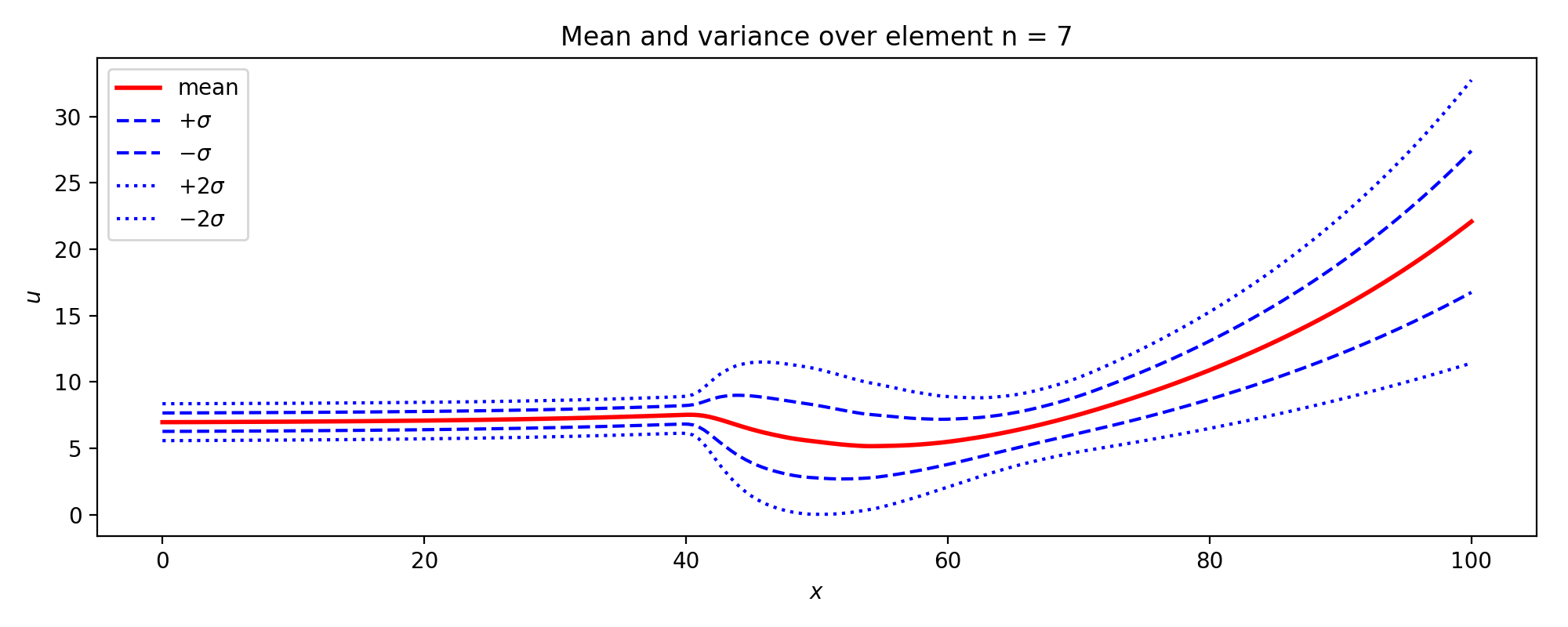}%
            {Mean and variance of the ROM solution at time $t=12$
              for $\cT^0$ (top) 
             and $\cT^7$ (bottom).The dashed and dotted lines
             indicate one and two standard deviations from the mean,
             respectively.}%
            {fig:u_mean_var}

\ezfiguretworow{width=0.9\textwidth}{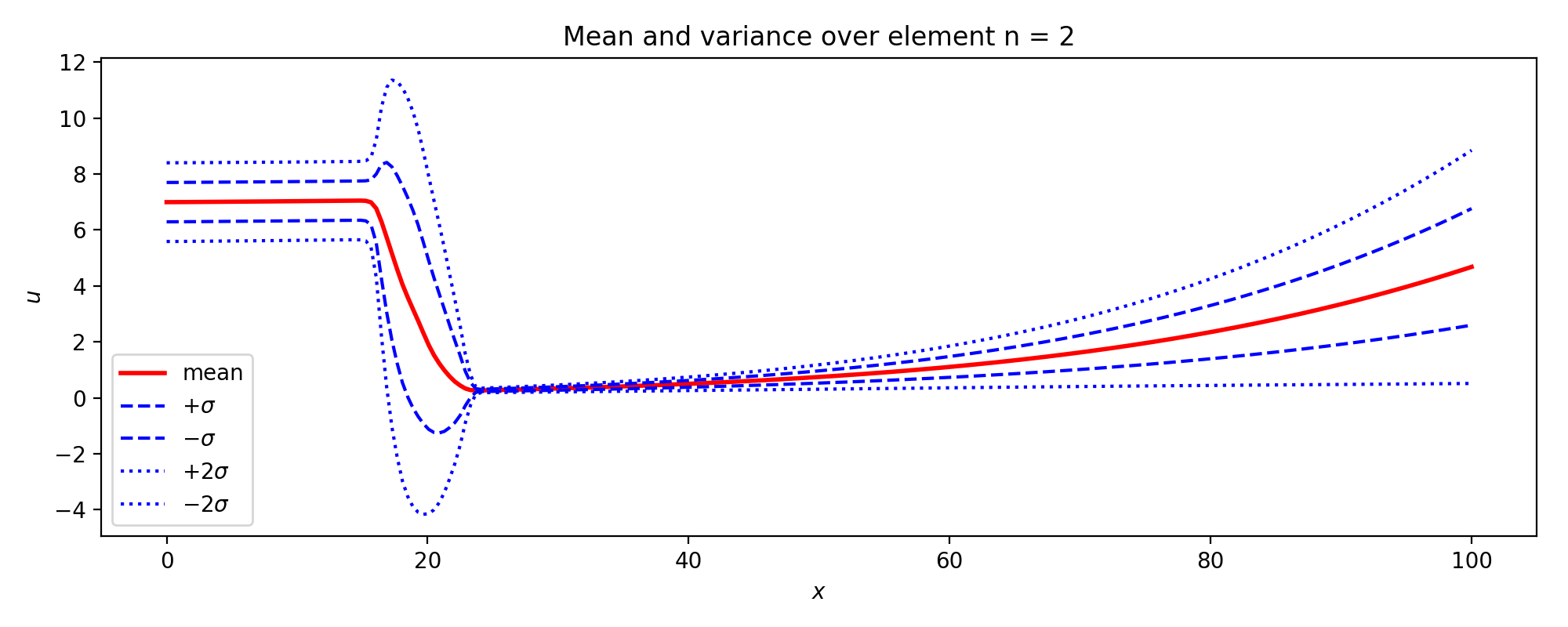}%
               {width=0.9\textwidth}{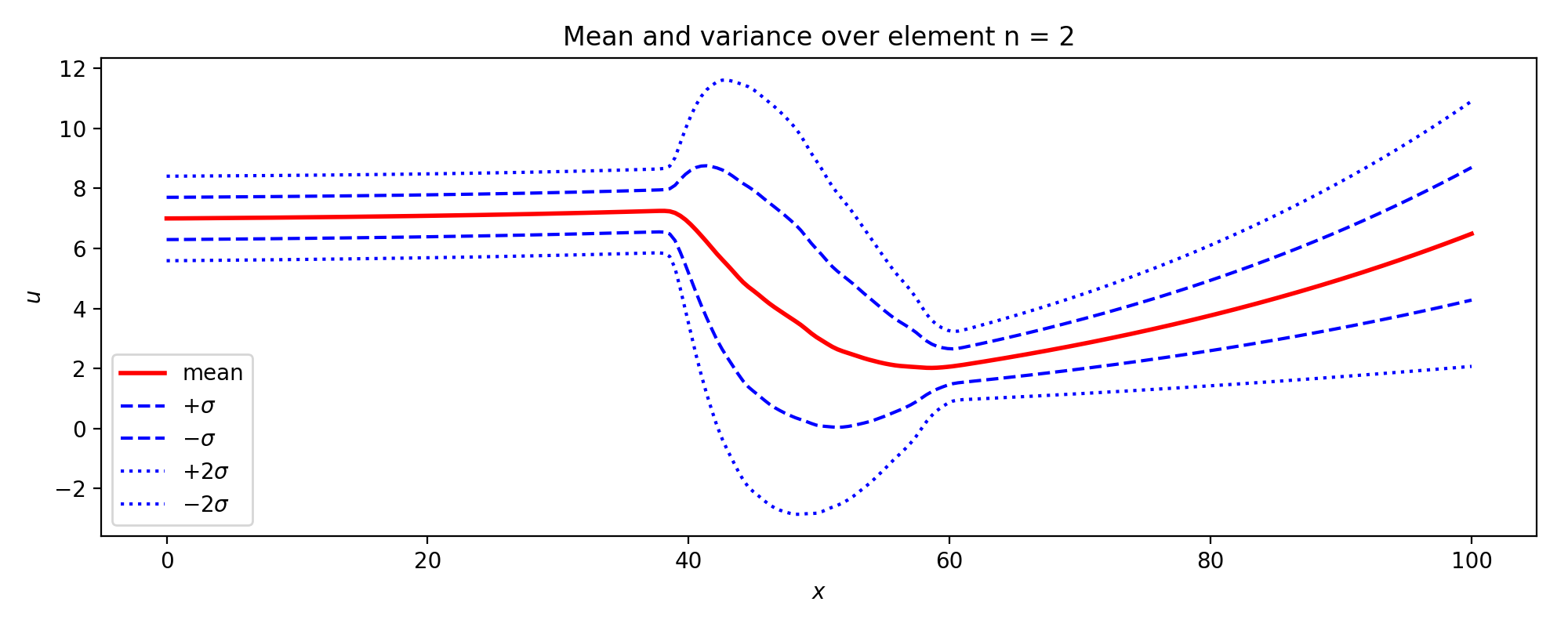}%
            {Mean and variance of the ROM solution at times $t=5$ (top) 
             and $t=12$ (bottom) for $\cT_2$. The dashed and dotted lines
             indicate one and two standard deviations from the mean,
             respectively.}%
            {fig:u_mean_var_time}

\item We compute the two QoIs, shock location and the shock height 
      at time $t = 12$, using the formulas for the centroid location of the 
      $n$-th piece $\mathfrak{q}^{-1}_u: \NN \to \mathfrak{c}(\cP(v))$ 
      \eqref{eq:ordering} and the integral of the $n$-th piece of $u$, 
      $\cP_2(n,u)$ \eqref{eq:P2},
      \beq
        \begin{aligned}
        \textrm{(shock location)}(t) &= \mathfrak{q}_u^{-1}(2), \\
        \textrm{(shock height)}(t) &= \int_0^{100} \cP_2(2,u) \dx.
        \end{aligned}
      \eeq
      A scatter plot for the two quantities, as functions of 
      the parameters $\mu_1$ and $\mu_2$ individually, are shown in 
      \cref{fig:scatter}.
\item The samples for the shock location and shock height can be 
      used to compute the kernel density estimate (KDE) between the two
      quantities, as shown in \cref{fig:joint_kde}.
      It is easy to see that while the correlation between the 
      two quantities are strong when both are small (smaller, slower shocks)
      the correlation weakens considerably when both are large 
      (large, faster shocks). 
\item Using linear regression, we can construct a surrogate surface
      for the shock height and shock location as functions of $\bmu$.
      Using $5$-th order polynomials, we can approximate the samples
      extremely well. This shows empirically that although the 
      solutions develop shocks and can be discontinuous, shock location
      can be a smooth function of the parameters.
      This observation is intimately related to the approximation properties
      we have discussed in \cref{sec:approx}, and has also been observed
      in \cite{li17}.
\end{itemize}

\ezfigure{width=1.00\textwidth}{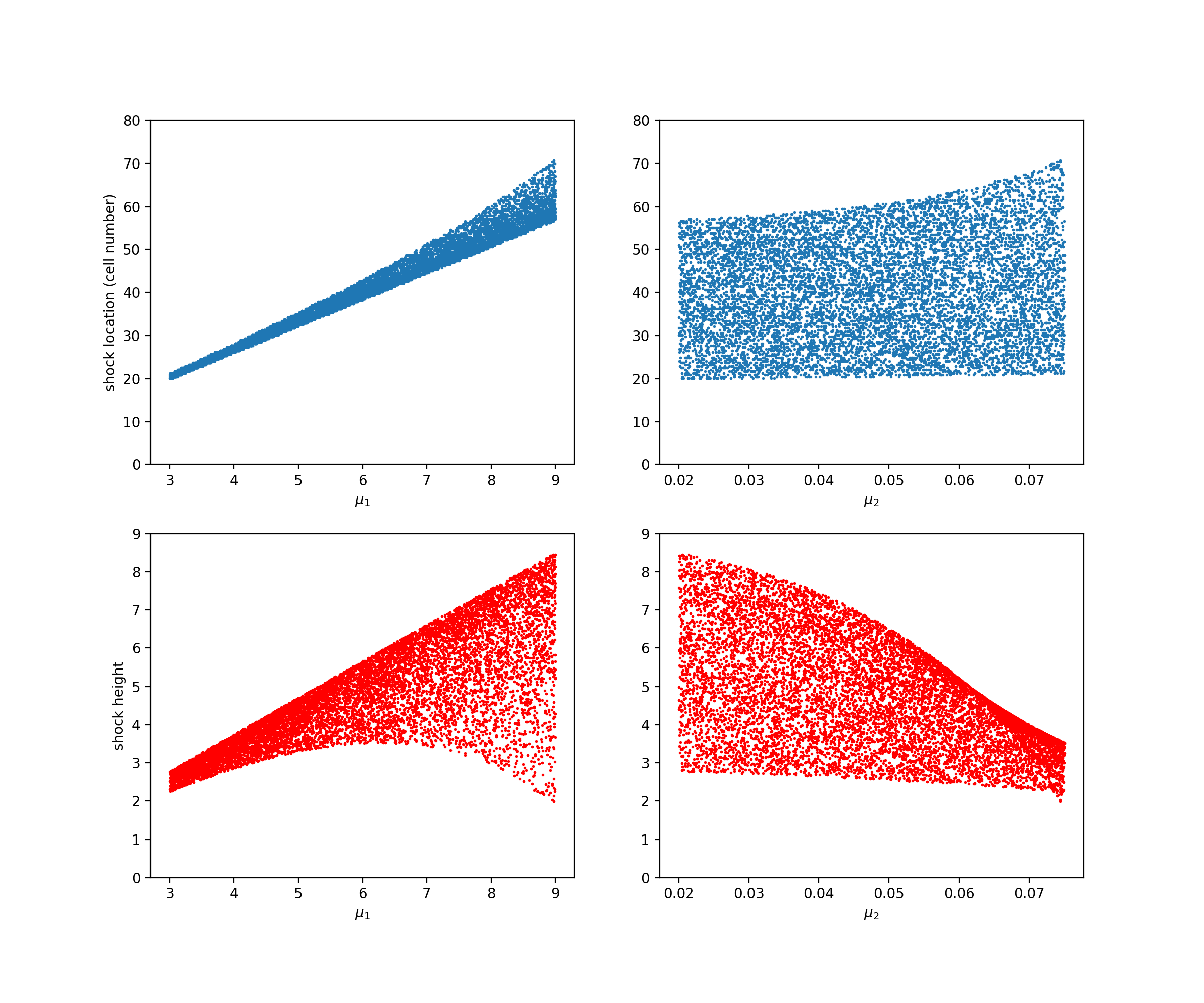}
         {Scatter plot of the shock height and shock location with respect
          to the parameters $\mu_1$ and $\mu_2$.}%
         {fig:scatter}

\ezfigure{width=0.9\textwidth}{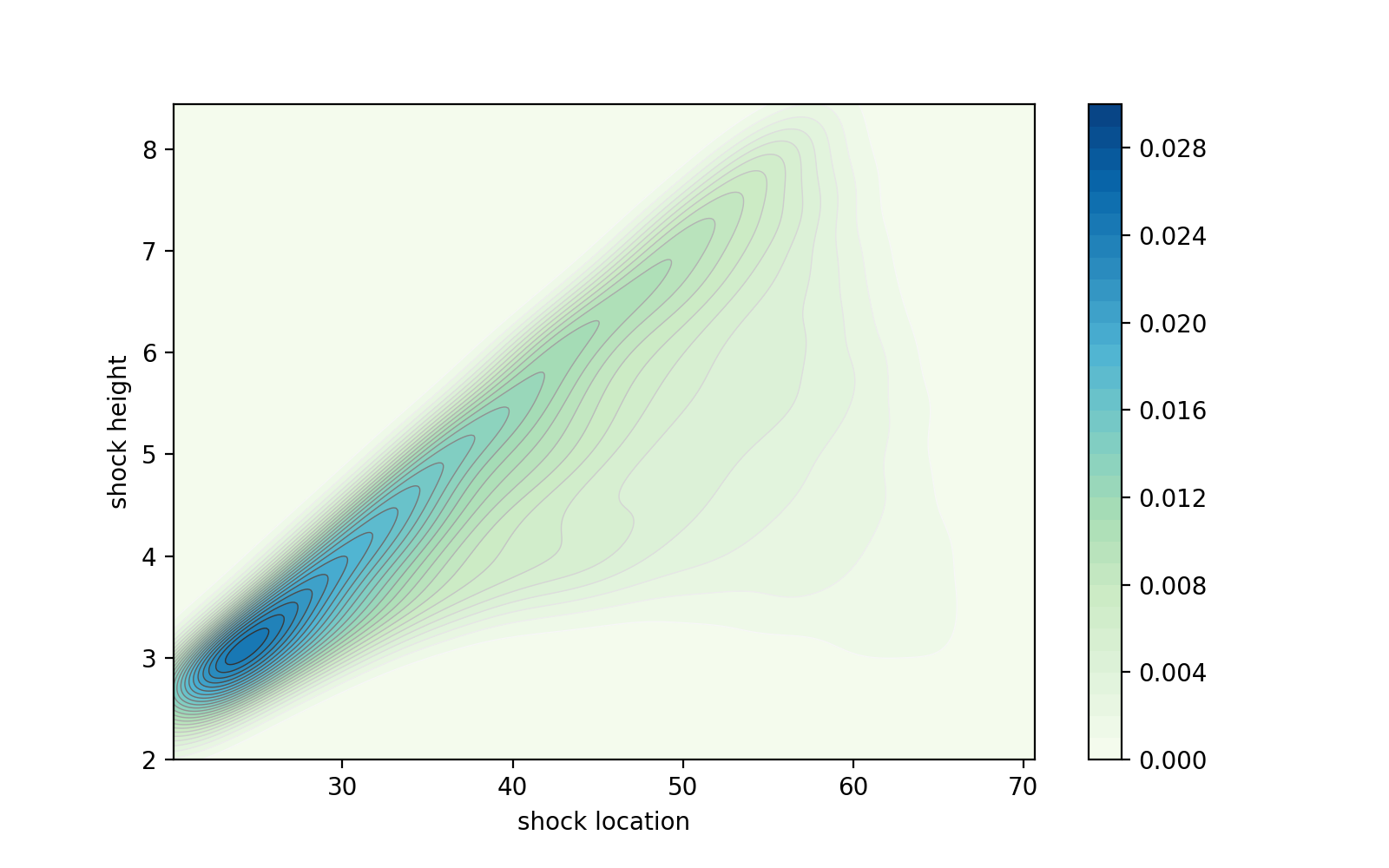}
         {Kernel density estimate (KDE) of the joint density between the 
          shock height and shock location.}%
         {fig:joint_kde}

\ezfiguretwo{width=0.45\textwidth}{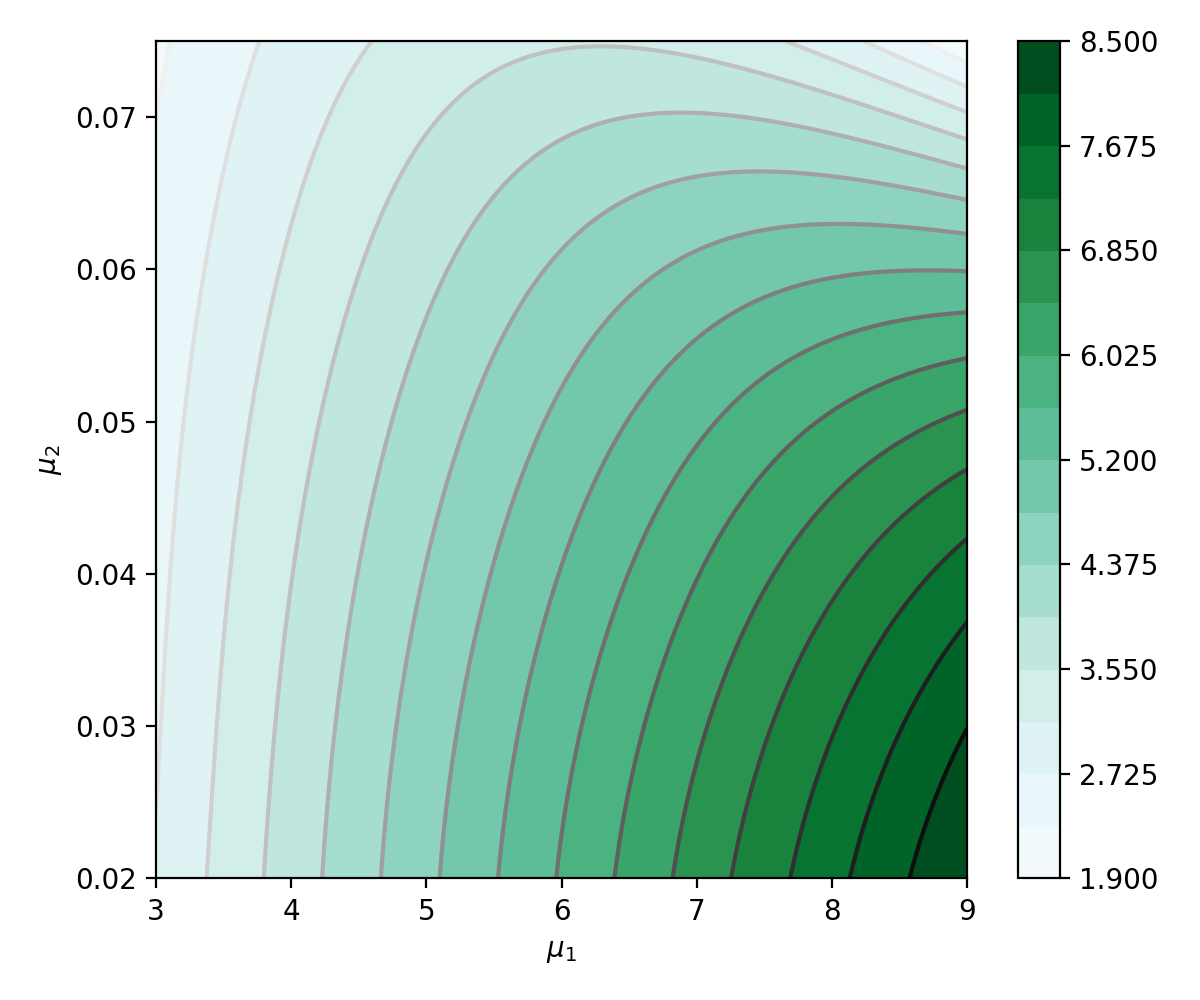}%
            {width=0.45\textwidth}{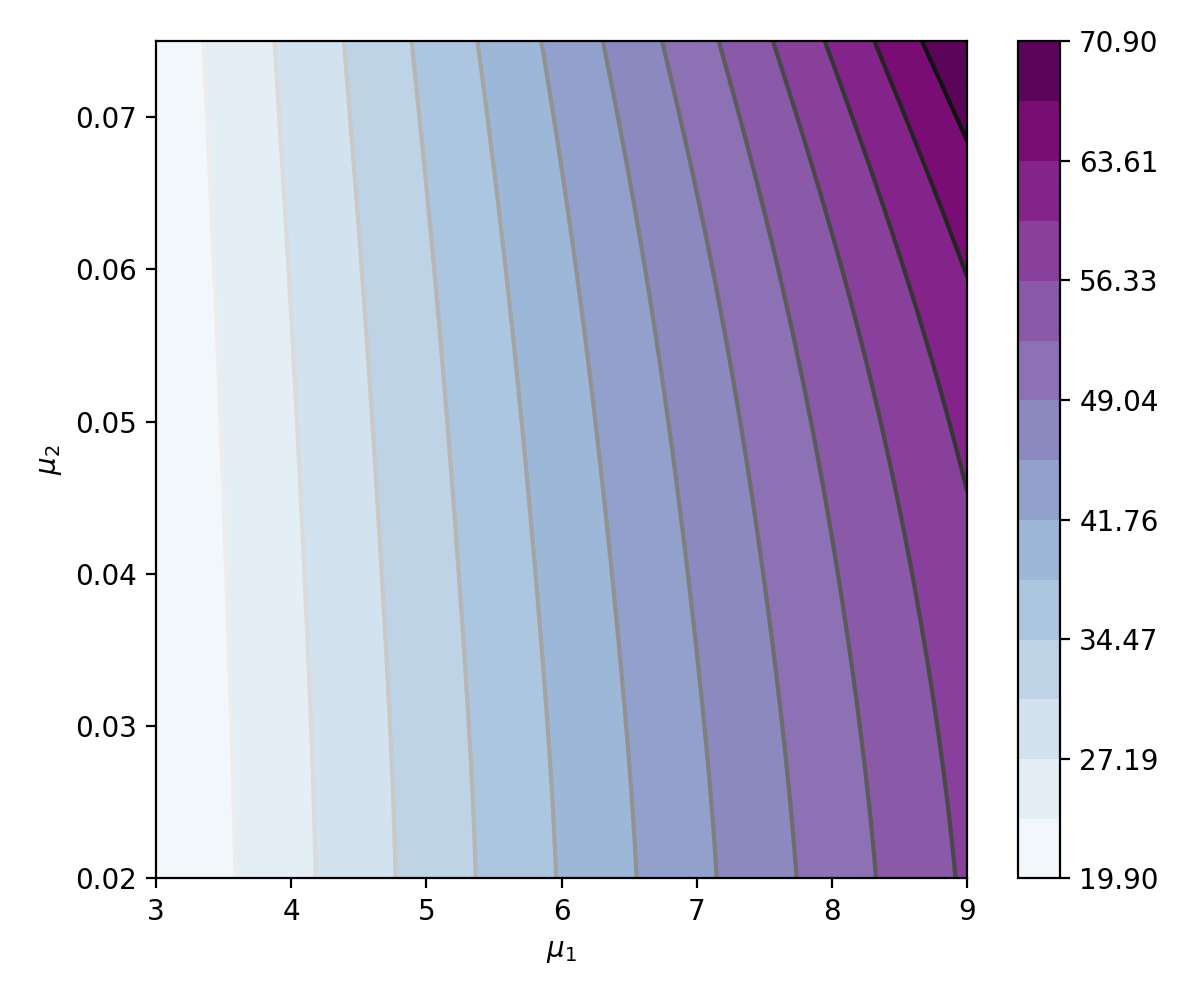}%
            {Surrogate surface for the shock height (left) and 
             shock location (right) at time $t=12$, 
             using polynomials of $\mu_1$ and $\mu_2$ orders up to $5$.}
            {fig:shock_hl}

\section{Conclusion and future work} \label{sec:conclusion}

We have proposed a technique for performing a model order reduction of a 
scalar hyperbolic conservation law (\cref{sec:def,sec:basis,sec:rom}),
successfully applied the technique  to a parametrized
Burgers' equation \eqref{eq:par_burg} and performed UQ tasks with 
the resulting ROM (\cref{sec:uq}). 
We have also discussed the approximation
properties of the displacement interpolation (\cref{sec:approx}) and proposed a
more general form for approximation \cref{eq:x_sum}.

We believe the themes that appear in this paper naturally leads to various 
topics for future investigation. 
The dynamics in terms of the solution-dependent variable in \cref{eq:dyn}
may be of use directly, lifting much of the burden in using the local basis
for the Galerkin update \cref{eq:rrom_update}.
The use of the triangulations \cref{eq:tri} to discretize the parameter space
is also not the sole option, and other approaches would be preferable if the
dimension of the parameter space is large. This would require suitably 
adapting the displacement interpolation \cref{eq:dip}. 

Application to examples that require entropy fixes, and
 application to systems of conservation laws, dealing with more difficult
source terms which appear for many applications for the shallow water equations 
\cite{berger11,leveque11}, generalization to multiple spatial dimensions, 
perhaps using the dimensional splitting approach using Radon transforms 
\cite{radonsplit} are topics that are directly related to this technique.
The implications of the signature condition
\cref{cond:signature} can be studied further, especially for multiple
dimensions. Further reduction for the time-stepping \cref{eq:rrom_update}
could be explored. Finally, development and analysis of approximations of the 
type \cref{eq:x_sum} and derivation of a more general approximation methods 
that can overcome the limitations of existing methods will also be pursued. 

\bibliographystyle{siamplain}

\end{document}

%% file: macros.tex
\newcommand{\ello}{{\ell}}

\newcommand{\cE}{\mathcal{E}}
\newcommand{\cF}{\mathcal{F}}

\newcommand{\cN}{\mathcal{N}}

\newcommand{\cO}{\mathcal{O}}
\newcommand{\cT}{\mathcal{T}}
\newcommand{\cM}{\mathcal{M}}
\newcommand{\cP}{\mathcal{P}}

\newcommand{\cS}{\mathcal{S}}
\newcommand{\cU}{\mathcal{U}}

\newcommand{\beq}{\begin{equation}}
\newcommand{\eeq}{\end{equation}}
\def\bals#1\eals{\begin{align*} #1 \end{align*}}
\def\bal#1\eal{\begin{align} #1 \end{align}}

\newcommand{\where}{\quad \text{ where } }
\newcommand{\Ffor}{\quad \text{ for }}

\newcommand\Dom\Omega

\newcommand\RR{\mathbb{R}}

\newcommand\NN{\mathbb{N}}

\newcommand\Lap\Delta

\newcommand\abs[1]{\left\lvert #1 \right\rvert}

\newcommand\dx{\,\mathrm{d}x}

\newcommand\dz{\,\mathrm{d}z}

\newcommand\ds{\,\mathrm{d}s}

\def\bpde#1\epde{\[\left\{\begin{aligned}#1\end{aligned}\right. \]}
\def\inbpde#1\inepde{\left\{\begin{aligned}#1\end{aligned}\right.}
\def\binpde#1\einpde{\left\{\begin{aligned}#1\end{aligned}\right.}

\def\cC{\mathcal{C}}

\def\cI{\mathcal{I}}

\def\cU{\mathcal{U}}

\def\cW{\mathcal{W}}

\def\cU{\mathcal{U}}

\def\p{\partial}

\def\bfW{\mathbf{W}}

\def\b0{\mathbf{0}}

\def\bmu{\mbox{\boldmath$\mu$}}

\def\bbmat{\begin{bmatrix}[r]}
\def\ebmat{\end{bmatrix}}

\def\balpha{\boldsymbol\alpha}

\newcommand{\barr}{\begin{array}}
\newcommand{\ea}{\end{array}}
\newcommand{\bea}{\begin{eqnarray}}
\newcommand{\eea}{\end{eqnarray}}
\newcommand{\bt}{\begin{table}}
\newcommand{\et}{\end{table}}

\DeclareMathOperator\Span{span}

\DeclareMathOperator\rg{range}

\DeclareMathOperator\sgn{sgn}

\theoremstyle{plain}
\theoremstyle{definition}

\newsiamthm{cond}{Condition}
\newsiamremark{remark}{Remark}

\numberwithin{equation}{section}

\newcommand\ezfigure[4]{%
\begin{figure}
\centering
\includegraphics[#1]{#2}
\caption{#3} \label{#4}
\end{figure}
}

\newcommand\ezfiguretwo[6]{%
\begin{figure}
\centering
\begin{tabular}{cc}
\includegraphics[#1]{#2}
&
\includegraphics[#3]{#4}
\end{tabular}
\caption{#5} \label{#6}
\end{figure}
}

\newcommand\ezfiguretworow[6]{%
\begin{figure}
\centering
\begin{tabular}{c}
\includegraphics[#1]{#2}
\\
\includegraphics[#3]{#4}
\end{tabular}
\caption{#5} \label{#6}
\end{figure}
}